\documentclass[11pt,reqno]{amsart}

\usepackage{amsmath,amsfonts,amsthm,amssymb,enumerate,color}
\usepackage{hyperref,graphicx}

\usepackage[a4paper,height=20cm,top=24mm,bottom=24mm,left=26mm,right=26mm]{geometry}

\newtheorem{thm}{Theorem}
\newtheorem{conj}{Conjecture}

\newcommand{\F}{{\mathbf F}}
\newcommand{\R}{{\mathbf R}}
\newcommand{\Z}{{\mathbf Z}}
\newcommand{\C}{{\mathbf C}}

\newcommand{\1}{{\bf 1} }

\renewcommand{\Re}{{\rm Re}}
\renewcommand{\Im}{{\rm Im}}
\renewcommand{\l}{\left}
\renewcommand{\r}{\right}

\newcommand{\rem}[1]{#1}
\newcommand{\remnum}[2]{#1}
\allowdisplaybreaks[1]

\begin{document}
\baselineskip 16pt
\parskip 8pt
\sloppy


\title{Euler Products beyond the Boundary}


\author[T.~Kimura]{Taro \textsc{Kimura}$^*$}
\thanks{$^*$Partially supported by JSPS Research Fellowships for Young
Scientists (Nos.~23-593, 25-4302)}

\author[S.~Koyama]{Shin-ya \textsc{Koyama}}

\author[N.~Kurokawa]{Nobushige \textsc{Kurokawa}}


\address{Mathematical Physics Laboratory, RIKEN Nishina Center, 2-1
Hirosawa, Wako, Saitama, 351-0198, Japan.}

\email{
\texttt{taro.kimura@riken.jp}}

\address{Department of Biomedical Engineering, Toyo University, 2100
Kujirai, Kawagoe, Saitama, 350-8585, Japan.}

\email{
\texttt{koyama@toyo.jp}}

\address{Department of Mathematics, Tokyo Institute of Technology,
2-12-1 Oh-okayama, Meguro-ku, Tokyo 152-8551, Japan.}

\email{
\texttt{kurokawa@math.titech.ac.jp}}



\begin{abstract}
\remnum{
We investigate the behavior of the Euler products of the Riemann zeta
 function and Dirichlet $L$-functions on the critical line.}{1}
A refined version of the Riemann hypothesis, which is named ``the Deep
 Riemann Hypothesis'' (DRH), is examined.
\remnum{
We also study various analogs for global function fields.
We give an interpretation for the nontrivial zeros from the viewpoint of
 statistical mechanics.}{1}
\end{abstract}





\maketitle

\section{Introduction}

Let $\chi$ be a primitive Dirichlet character with conductor $N$.
The Dirichlet $L$-function is expressed by \remnum{an}{2} Euler product
\begin{equation}\label{1.1}
L(s,\chi)=\prod_p(1-\chi(p)p^{-s})^{-1},
\end{equation}
where $p$ runs through all primes.
The product \eqref{1.1} is absolutely convergent \remnum{for}{3} $\Re(s)>1$.
It is known that $L(s,\chi)$ has a meromorphic continuation to all $s\in\C$,
which is entire if $\chi\ne\1$, and has a simple pole at $s=1$ if $\chi=\1$.

In this paper we study the values $L(s,\chi)$ beyond the boundary
$\Re(s)=1$ of the absolute convergence region $\Re(s)>1$
from the viewpoint of its relation to the values of the Euler product.
Few results are known \remnum{in}{4} this context. The classical results 
concerning the fact that
the Euler product \eqref{1.1} converges to $L(1+it,\chi)$ $(t\in\R,\ t\ne0)$
\remnum{
can be found in textbooks for either}{5}
$\chi=\1$ (\cite{T} Chapter 3) or $\chi\ne\1$ (\cite{M}).
The only work we \remnum{could}{6} find beyond this
is that of Goldfeld \cite{G}, Kuo-Murty \cite{KM} and Conrad \cite{C}.
Goldfeld \cite{G} and Kuo-Murty \cite{KM} dealt with
the $L$-functions of elliptic curves at $s=1$\remnum{,}{7}
with their results supporting the Birch and Swinnerton-Dyer conjecture.
Conrad \cite{C} treated more general Euler products for $\Re(s)\ge1/2$.

The (generalized) Riemann Hypothesis (GRH) for $L(s,\chi)$ asserts that
$L(s,\chi)\ne 0$ in $\Re(s)>1/2$. When $\chi\ne\1$, it is equivalent to
the following conjecture \remnum{\cite{C}}{8}.

\begin{conj}\label{conj1}
If $\chi\ne\1$, then for $\Re(s)>1/2$ we have
$$
L(s,\chi)=\lim_{n\to\infty}\prod_{p\le n}(1-\chi(p)p^{-s})^{-1},
$$
where the product is taken over all primes $p$ satisfying $p\le n$.
\end{conj}

Note that the order of primes which participate in the product is important, because
it is not absolutely convergent.


\begin{conj}[Deep Riemann Hypothesis (DRH)]\label{conj2}
If $\chi\ne\1$ and $L(s,\chi)\ne0$ with $\Re(s)=\frac12$, we have
\[
 \lim_{n\to\infty}\prod_{p\le n}(1-\chi(p)p^{-s})^{-1}
=L(s,\chi)\times
\begin{cases}
\sqrt 2 & (s=\frac12{\text{ and }}\chi^2=\1)\\
1& (\text{otherwise})
\end{cases},
\]
where the product is taken over all primes $p$ satisfying $p\le n$.
\end{conj}

We call Conjecture~\ref{conj2} the Deep Riemann Hypothesis,
a deeper modification of Conjecture~\ref{conj1}, literally because we reach
the boundary of the domain $\Re(s)>1/2$ given in Conjecture~\ref{conj1},
and logically because Conjecture~\ref{conj2} implies Conjecture~\ref{conj1}.
Indeed, if we denote
$$
\psi(x, \chi)=\sum_{m=1}^\infty \sum_{p:\ p^m\le x}\chi(p)\log p,
$$
Conjecture~\ref{conj1} is equivalent to
$$
\psi(x,\chi)=O(\sqrt x (\log x)^2),
$$
while Conjecture~\ref{conj2} is equivalent to
$$
\psi(x,\chi)=o(\sqrt x \log x)
$$
\remnum{by Conrad \cite{C} Theorem 6.2.}{9}

The prototype version of this Conjecture~\ref{conj2} was proposed in~\cite{C}.
For a generalization of Conjecture \ref{conj2} to the case including
$\chi=\1$, see Akatsuka~\cite{A}.

It is an easy task to obtain numerical support of Conjecture \ref{conj2},
since the convergence of the left hand side is fairly fast.

This kind of process, introducing a parameter to define a finite analogue
and then taking it to infinity, is often used in physics when it is
difficult to analyze the infinite system directly.
One can investigate how to approach infinity by analyzing the deviation
from the result in the desirable limit.
For example, in order to study \remnum{the}{10} asymptotic behavior in
\remnum{an}{10} infinite volume system, it is convenient to introduce a
system of some finite size $\Lambda$, and then estimate a correction by
analyzing a differential equation in terms of $\Lambda$, which is
\remnum{the}{10} so-called {\em renormalization group equation}.

The situation for the Riemann zeta and the Dirichlet $L$-functions seems
quite similar: the difficulty \remnum{with}{10} these functions
\remnum{lies}{10} essentially involved \remnum{in treating
infinity}{10}, so that convergency of the Euler product is nontrivial.
In this paper we numerically examine {\em the finite-size corrections}
to the zeta and $L$-functions appearing in the finite analog, based on
the analogy between nontrivial zeros and eigenvalues of a certain
infinite dimensional matrix or critical phenomena observed around a
phase transition point.

\section{Function Field Analogs}

In this section, we prove an analog of Conjecture \ref{conj2} for
function fields of one variable over a finite field.
\remnum{
The theory of zeta and $L$-functions over such function fields are seen,
for example, in the textbook of Rosen \cite{Rosen2}.}{12}

Let $\F_q$ be the finite field of $q$ elements.
We fix a conductor $f(T)\in\F_q[T]$ and introduce a ``Dirichlet'' character
$$\chi:\ (\F_q[T]/(f))^\times \to \C^\times,$$
\remnum{
which is extended to $\F_q[T]$ by $\chi (h)=0$ for $h$ such that
$(h,f)\ne (1)$.}{13}
We define the ``Dirichlet'' $L$-function by the Euler product:
$$L_{\F_q(T)}(s,\chi)=\prod_h (1-\chi(h)N(h)^{-s})^{-1},$$
where $h=h(T)\in\F_q[T]$ runs through monic irreducible polynomials,
and $N(h)=q^{\deg h}$.
\remnum{In}{14} the celebrated work of Kornblum \cite{K}, it is proved that 
the above Euler product is absolutely convergent in $\Re(s)>1$,
and is a polynomial in $q^{-s}$ of degree less than 
\remnum{$\deg (f)-2$ if $\chi\ne\1$~\cite{Weil2}.}{15}

We prove the following theorem.
\begin{thm}[DRH over function fields]\label{thm1}
Let $q$, $f$ and $\chi$ be as above.
Put $K=\F_q(T)$ and assume $\chi\ne\1$.
Then the following (1) and (2) are true.
\begin{enumerate}[\rm(1)]
\item For $\Re(s)>1/2$, we have
$$\lim_{n\to\infty}\prod_{\deg h\le n}(1-\chi(h)N(h)^{-s})^{-1}=L_K(s,\chi).$$
\item For $t\in\R$ with $L_K(\frac12+it,\chi)\ne0$, it holds that
\[
\lim_{n\to\infty}\prod_{\deg h\le n}(1-\chi(h)N(h)^{-\frac12-it})^{-1}
=L_K\l(\frac12+it,\chi\r)\times
\begin{cases}\sqrt 2 & (\chi^2=\1,\ t\in\frac\pi{\log p}\Z)\\1& (\text{otherwise})\end{cases}.
\]
\end{enumerate}
\end{thm}

\medskip
{\it Proof of Theorem \ref{thm1}.}~
\remnum{We prove (2) first.}{16}
We estimate the product
$$E_n=\prod_{\deg h\le n}\l(1-\chi(h)N(h)^{-\frac12-it}\r)^{-1}$$
by dealing with its logarithm
$$
\log E_n=\sum_{\deg h\le n}\sum_{k=1}^\infty \frac{\chi(h)^k}{k}q^{-k(\frac12+it)\deg h}.
$$
We divide the sum into three parts as $$\log E_n=A(n)+B(n)+C(n)$$ with
\begin{align*}
A(n)&=\sum_{k=1}^\infty \sum_{\deg h\le n/k}\frac{\chi(h)^k}{k}q^{-k(\frac12+it)\deg h},\\
B(n)&=\sum_{n/2\le \deg h\le n}\frac{\chi(h)^2}{2}q^{-2(\frac12+it)\deg h},\\
C(n)&=\sum_{k=3}^\infty \sum_{n/k< \deg h\le n}\frac{\chi(h)^k}{k}q^{-k(\frac12+it)\deg h}.
\end{align*}

By the above mentioned Kornblum's theorem, we put
$$
L_K(s,\chi)=\prod_{j=1}^r (1-\lambda_jq^{-s})
$$
with $|\lambda_j|=\sqrt q$ or 1~\remnum{\cite{Deligne}\cite{Grothendieck}\cite{Weil1}.}{17} 
Then by taking the logarithmic derivatives of
$$
\prod_h (1-\chi(h)N(h)^{-s})^{-1}=\prod_{j=1}^r (1-\lambda_jq^{-s}) \quad (\Re(s)>1)
$$
and comparing the coefficients of $q^{-sk}$, we have
$$
\sum_{(\deg h)|k}(\deg h)\chi(h)^{\frac k{\deg h}}
=-\sum_{j=1}^r \lambda_j^k\quad (k\ge 1).
$$
By this identity, the first partial sum $A(n)$ is calculated as
\begin{align*}
A(n)
&=\sum_{k\le n}\frac{q^{-(\frac12+it)k}}k\sum_{(\deg h)|k}(\deg h)\chi(h)^{\frac k{\deg h}}\\
&=-\sum_{j=1}^r \sum_{k=1}^n \frac1k\l(\frac{\lambda_j}{q^{\frac12+it}}\r)^k.
\end{align*}
By the \remnum{Deligne's}{18} theorem we have
$\l|\frac{\lambda_j}{q^{\frac12+it}}\r|\le1$ and the assumption
$L_K(\frac12+it,\chi)\ne0$ tells that
$\frac{\lambda_j}{q^{\frac12+it}}\ne1$.
\remnum{
Then by the Taylor expansion for $\log (1-x)$,}{19} it holds that
\begin{eqnarray}
 \lim_{n\to\infty}A(n)
& = & \sum_{j=1}^r \log \l(1-\frac{\lambda_j}{q^{\frac12+it}}\r)
\nonumber\\ 
& = & \log L_K\l(\frac12+it,\chi\r).
 \nonumber
\end{eqnarray}

Next for estimating $B(n)$, we use the 
\remnum{generalized Mertens' theorem~\cite{Rosen1}}{20} that
$$\sum_{\deg h<n}\frac1{N(h)}\sim \log n
\qquad (n\to\infty).$$
When $\chi^2=\1$ and $t\in\frac\pi{\log q}\Z$, we compute that
\begin{align*}
B(n)
&=\frac12\sum_{n/2\le \deg h\le n} q^{-(1+2it)\deg h}\\
&=\frac12\l(\sum_{1\le \deg h\le n} q^{-(1+2it)\deg h}\r.
\l.-\sum_{1\le \deg h< n/2} q^{-(1+2it)\deg h}\r)\\
&=\frac12\l(\l(\log n+C+O(n^{-1})\r)
-\l(\log \frac n2+C+O(n^{-1})\r)\r)\\
&=\frac12\l(\log 2+O(n^{-1})\r).
\end{align*}
Hence 
$$\lim_{n\to\infty}B(n)=\log \sqrt 2.$$
\remnum{
In all other cases it holds that $B(n)\to 0$ as $n \to \infty$.}{21}

Finally, $C(n)\to 0$ as $n\to\infty$
\remnum{by a similar argument to Lemma 3.1 in \cite{C}.}{22}

\remnum{
For proving (1), we use the decomposition into $A(n)$ and $B(n)+C(n)$, in
place of that into $A(n)$, $B(n)$ and $C(n)$ above.
In this case both $B(n)$ and $C(n)$ are concerning absolutely convergent
series like $C(n)$ in the proof of (2).
Thus $B(n)+C(n) \to 0$ as $n \to \infty$.}{16}
\hfill\qed

Conjecture \ref{conj2} and Theorem \ref{thm1} are generalized to
automorphic $L$-functions \remnum{by}{23} Lownes \cite{L}.

The following theorems are for the case of the trivial character.
\begin{thm}\label{thm2}
Let $X$ be a projective smooth curve over $\F_q$. Then
$$
\lim_{n\to \infty}
\prod_{N(x)\le q^n}(1-N(x)^{-1/2})^{-1} \cdot 
\rem{\exp\l(
-\sum_{l=1}^n\frac{q^{l/2}}l\r)}
 = \sqrt2 \l(\sqrt q-1\r)
 \l|\rem{\zeta\l(X,\frac12\r)}
 \r| .
$$
\end{thm}

Notice that 
$$
\sum_{l=1}^n\frac{q^{l/2}}l
=\frac{\log q}{1-q^{-1}}\int_1^{q^n}\frac{d_q(u)}{\sqrt u\log u},
$$
where
$$
\int_1^{q^n}f(u) d_q(u)
=\sum_{l=1}^n f(q^l)(q^l-q^{l-1})
$$
is \remnum{Jackson's $q$-integral~\cite{KC}\cite{Jackson}.}{25} 
Thus, it is considered as a
``modified $q$-logarithmic integral.'' 
The situation is 
\remnum{extended}{25} to the case of the Riemann zeta function studied by
Akatsuka~\cite{A},
where a ``modified logarithmic integral'' appears.

\medskip
{\it Proof of Theorem \ref{thm2}.}~
Let $g$ be the genus of the curve $X$.
\remnum{
By Deligne's theorem \cite{Deligne}}{26} there exist
$\alpha_j\in\C$ with $|\alpha_j|=\sqrt q$ for $j=1,2,3,...,g$
such that 
$$
\zeta(X,s)=\frac{\prod_{j=1}^g
(1-\alpha_jq^{-s})(1-\overline{\alpha_j}q^{-s})}
{(1-q^{-s})(1-q^{1-s})}.
$$
Note that $\alpha_j\ne\sqrt q$, because
$\alpha_j+\overline{\alpha_j}\in\Z$.
Thus we have
$$
\zeta\l(X,\frac12\r)=\frac{\prod_{j=1}^g
(1-\alpha_jq^{-1/2})(1-\overline{\alpha_j}q^{-1/2})}
{(1-q^{-1/2})(1-q^{1/2})}.
$$
On the other hand we compute
\begin{align*}
\lefteqn{\log\l(\prod_{N(x)\le q^n}(1-N(x)^{-1/2})^{-1}\r)}\\
&=\log\prod_{\deg(x)\le n}\l(1-q^{-\frac{\deg(x)}2}\r)^{-1}\\
&=\sum_{\deg(x)\le n}\sum_{k=1}^\infty\frac{q^{-\frac{k\deg(x)}2}}k\\
&=\sum_{\substack{\rem{k,n}
\\k\deg(x)\le n}}\frac{q^{-\frac{k\deg(x)}2}}k
+\frac12\sum_{\frac n2<\deg(x)\le n}q^{-\deg (x)}
+\sum_{k=3}^\infty\frac1k\sum_{\frac nk<\deg(x)\le n}q^{-\frac{k\deg(x)}2}.
\end{align*}
When $n\to\infty$, the second term tends to $\frac12\log2$ by the
generalized Mertens' theorem~\cite{Rosen1}, 
and the third term goes to 0, because we have
$\sum_{x\in|X|}N(x)^{-\alpha}<\infty$ for any $\alpha>1$.
The first term is calculated as follows.
\begin{align*}
\sum_{\substack{\rem{k,n}
\\k\deg(x)\le n}}\frac{q^{-\frac{k\deg(x)}2}}k
&=\sum_{l=1}^n\frac1l\l(\sum_{\deg(x)|l}\deg(x)\r)q^{-l/2}\\
&=\sum_{l=1}^n\frac{|X(\F_{q^l})|}lq^{-l/2}\\
&=\sum_{l=1}^n\frac{q^l+1-\sum_{j=1}^g(\alpha_j^l+\overline{\alpha_j^l})}lq^{-l/2}\\
&=\sum_{l=1}^n\frac{q^{l/2}}l+\sum_{l=1}^\infty\frac{1-\sum_{j=1}^g(\alpha_j^l+\overline{\alpha_j^l})}lq^{-l/2}+o(1)\\
&=\sum_{l=1}^n\frac{q^{l/2}}l+\log\frac{\prod_{j=1}^g(1-\alpha_jq^{-1/2})(1-\overline{\alpha_j}q^{-1/2})}
{1-q^{-1/2}}+o(1),
\end{align*}
\remnum{
where we used 
the fact that $|\alpha_j|=\sqrt q$  $(\alpha_j\ne\sqrt q)$
for convergence of the Taylor expansion of the logarithms.}{27}
Therefore it holds that
\begin{multline*}
\log\l(\prod_{N(x)\le q^n}(1-N(x)^{-1/2})^{-1}\r)\\
=\sum_{l=1}^n\frac{q^{l/2}}l+\frac12\log2+
\log\frac{\prod_{j=1}^g(1-\alpha_jq^{-1/2})(1-\overline{\alpha_j}q^{-1/2})}{1-q^{-1/2}}+o(1).
\end{multline*}
Hence
$$
\prod_{N(x)\le q^n}(1-N(x)^{-1/2})^{-1}
\sim\exp\l(\sum_{l=1}^n\frac{q^{l/2}}l\r)\sqrt2\l(\sqrt q-1\r)
\l|\zeta\l(X,\frac12\r)\r|.
$$
\hfill\qed

\remnum{
Theorem \ref{thm2} is the ``deeper analogue'' for smooth curves of the
following Theorem \ref{thm3} for proper smooth schemes, which in its
turn is a function field analogue of Mertens' theorem \cite{Rosen1}.}{30}
In the situation of Theorem \ref{thm2}, it holds that
$$
\prod_{N(x)\le t}(1-N(x)^{-1})^{-1}
\sim\l(\mathrm{Res}_{s=1}\zeta(X,s)\r)e^\gamma \log t
$$
as $t\to\infty$.

\begin{thm}\label{thm3}
Let $X$ be a proper smooth scheme over $\F_p$. 
Then we have
$$
\prod_{N(x)\le t}(1-N(x)^{-\dim(X)})^{-1}
\sim\l(\mathrm{Res}_{s=\dim(X)}\zeta(X,s)\r)e^\gamma\log t
$$
as $t\to\infty$.
\end{thm}

{\it Proof of Theorem \ref{thm3}.}~
\begin{align*}
\lefteqn{\log\l(\prod_{N(x)\le q^n}(1-N(x)^{-\dim(X)})^{-1}\r)}\\
&=\log\l(\prod_{\rem{\deg(x)}
 \le n}(1-q^{-\dim(X)\deg(x)})^{-1}\r)\\
&=\sum_{\deg(x)\le n}\sum_{k=1}^\infty\frac{q^{-\dim(X)k\deg(x)}}k\\
&=\sum_{\substack{\rem{k,n}
 \\k\deg(x)\le n}}\frac{q^{-\dim(X)k\deg(x)}}k
+\sum_{k=2}^\infty\frac1k\sum_{\frac nk<\deg(x)\le n} q^{-\dim(X)k\deg(x)}.
\end{align*}
The second term goes to 0 as $n\to\infty$, because we have
$\sum_{x\in|X|}N(x)^{-\alpha}<\infty$ for any $\alpha>\dim(X)$.
The first term is calculated as follows.
By putting $l=k\deg(x)$, we compute
\begin{align}\label{GD}
\sum_{k\deg(x)\le n}\frac{q^{-\dim(X)k\deg(x)}}k
&=\sum_{l=1}^n\frac1l\l(\sum_{\deg(x)|l}\deg(x)\r)q^{-\dim(X)l}.
\end{align}
By the results of \remnum{Grothendieck \cite{Grothendieck} and Deligne
\cite{Deligne}}{34}, there exist
$\alpha_i$, $\beta_j\in\C$ with
$|\alpha_i|$, $|\beta_j|<q^{\dim(X)}$ such that
\begin{align*}
\sum_{\deg(x)|l}\deg(x)
&=|X(\F_{q^l})|\\
&=q^{l\dim(X)}+\sum_j\beta_j^l-\sum_i\alpha_i^l.
\end{align*}
Hence
\begin{align*}
\eqref{GD}
&=\sum_{l=1}^n\frac1l+\sum_{l=1}^n\frac1l
\l(\sum_{j}\l(\frac{\beta_j}{q^{\dim(X)}}\r)^l-\sum_{i}\l(\frac{\alpha_i}{q^{\dim(X)}}\r)^l\r)\\
&=\log n+\gamma+\log
\frac{\prod_i(1-\alpha_i q^{-\dim(X)})}{\prod_j(1-\beta_j q^{-\dim(X)})}+o(1),
\end{align*}
as $n\to\infty$. 
Since
$$
\zeta(X,s)
=\frac{\prod_i(1-\alpha_iq^{-s})}{(1-q^{\dim(X)-s})\prod_j(1-\beta_j q^{-s})},
$$
we see that $s=\dim(X)$ is the largest pole of $\zeta(X,s)$, which is simple with
$$
\mathrm{Res}_{s=\dim(X)}\zeta(X,s)=\frac1{\log q}\cdot
\frac{\prod_i(1-\alpha_i q^{-\dim(X)})}{\prod_j(1-\beta_j q^{-\dim(X)})}.
$$
Taking all terms into account, we conclude that
\begin{align*}
\prod_{N(x)\le q^n}(1-N(x)^{-\dim(X)})^{-1}
&\sim ne^\gamma(\log q)\mathrm{Res}_{s=\dim(X)}\zeta(X,s)\\
&=(\log q^n)\cdot e^\gamma \cdot\mathrm{Res}_{s=\dim(X)}\zeta(X,s).
\end{align*}
\hfill\qed

We conjecture that Theorem \ref{thm3} would hold for general schemes:
\begin{conj}
Let $X$ be a proper smooth scheme over $\Z$. Then
$$
\prod_{N(x)\le t}(1-N(x)^{-\dim(X)})^{-1}
\sim\l(\mathrm{Res}_{s=\dim(X)}\zeta(X,s)\r)e^\gamma\log t
$$
as $t\to\infty$.
\end{conj}

\section{Numerical Calculations}
In this section we show some \remnum{numerical data supporting}{35,36}
the Deep Riemann Hypothesis (Conjecture \ref{conj2}).
\remnum{If this conjecture is true,}{a} the partial Euler product
$$
L_x(s,\chi)=\prod_{p\le x}(1-\chi(p)p^{-s})^{-1},
$$
converges to $L(s,\chi)$ or $\sqrt2 L(s,\chi)$ as $x\to\infty$ even
on the critical line $\Re(s)=1/2$.
We formally put $L_x(s,\chi)=L(s,\chi)$ for $x=\infty$.

First we give Table \ref{Table1}, which shows the accuracy of Conjecture \ref{conj2} at $s=1/2$.
We find that the ratio of $\sqrt2L(\frac12,\chi)$ and
$L_x(\frac12,\chi)$ is almost equal to 1 for $x=10^7$,
when $\chi$ is quadratic.

\begin{table}[h]
 \begin{center}
\begin{tabular}{|c|c|c|c|}\hline
$d$
&
$\sqrt{2}L$
& 
$E$
& $(\sqrt{2}L)/E$
\\ \hline\hline
$-3$  & $0.680049$ & $0.688002$ & $0.988440$ \\
$-4$  & $0.944258$ & $0.945909$ & $0.998254$ \\
$5$   & $0.327745$ & $0.320619$ & $1.022223$ \\ 
$-7$  & $1.621517$ & $1.640320$ & $0.988536$ \\
$8$   & $0.528479$ & $0.539992$ & $0.978680$ \\ 
$-8$  & $1.556230$ & $1.521663$ & $1.022716$ \\
$-11$ & $1.402301$ & $1.342967$ & $1.044181$ \\
$12$  & $0.705066$ & $0.729170$ & $0.966942$ \\ 
$13$  & $0.621678$ & $0.618558$ & $1.005044$ \\ 
$-15$ & $2.612093$ & $2.791265$ & $0.935809$ \\
$17$  & $1.020601$ & $1.066235$ & $0.957201$ \\ 
$-19$ & $1.137621$ & $1.173052$ & $0.969795$ \\
$-20$ & $2.375413$ & $2.356696$ & $1.007942$ \\
$21$  & $0.703235$ & $0.724051$ & $0.971250$ \\ 
$-23$ & $3.472406$ & $3.320551$ & $1.045732$ \\
$24$  & $1.003325$ & $1.057376$ & $0.948881$ \\ 
$-24$ & $2.223023$ & $2.130498$ & $1.043428$ \\
$28$  & $1.162994$ & $1.199957$ & $0.969196$ \\ 
$29$  & $0.658655$ & $0.683281$ & $0.963958$ \\ 
\hline
\end{tabular}
\end{center}
\caption{$L:=L\left(\frac{1}{2},\left(\frac{d}{\cdot}\right)\right)$,
 $E:=\displaystyle\prod_{p\leq
 10^7}\left(1-\left(\frac{d}{p}\right)\frac{1}{\sqrt{p}}\right)^{-1}.$}
 \label{Table1}
\end{table}

In what follows we put
$\chi_{7a}$ and $\chi_{7b}$
to be the character $\chi$ modulo 7 
with $\chi^2\ne\1$ and $\chi^2=\1$, respectively.
Namely, if we define the character $\chi$ modulo 7 by giving the value
at the primitive root $3\in\Z/7\Z$, we define
$\chi_{7a}(3)=\exp(\pi\sqrt{-1}/3)$ and $\chi_{7b}(3)=-1$.
We also denote by $\chi_3$ the nontrivial character modulo 3,
which satisfies $\chi_3^2=1$.

\begin{figure}[h]
 \begin{center}
  \includegraphics[height=10em]{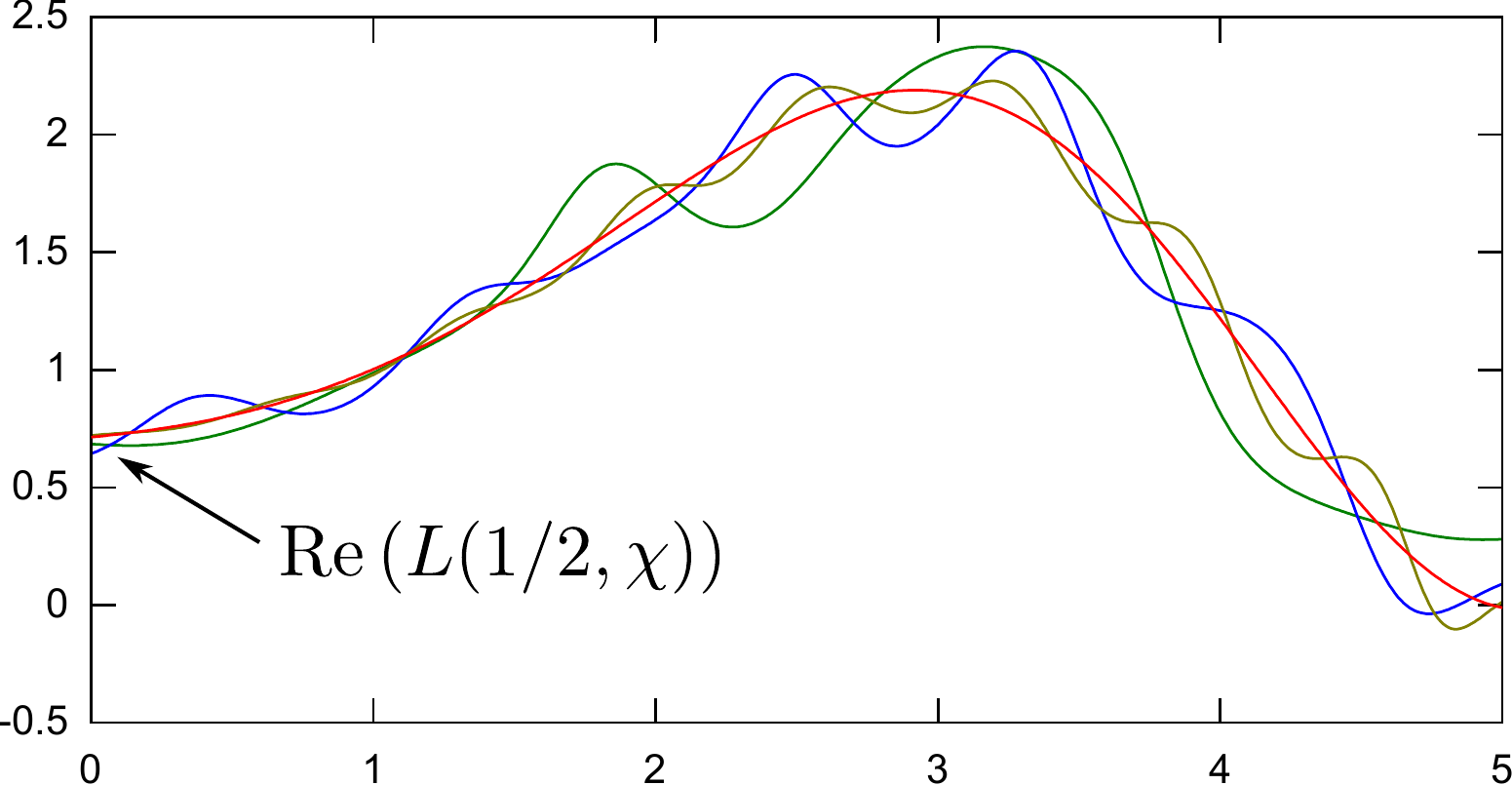} \qquad
  \includegraphics[height=10em]{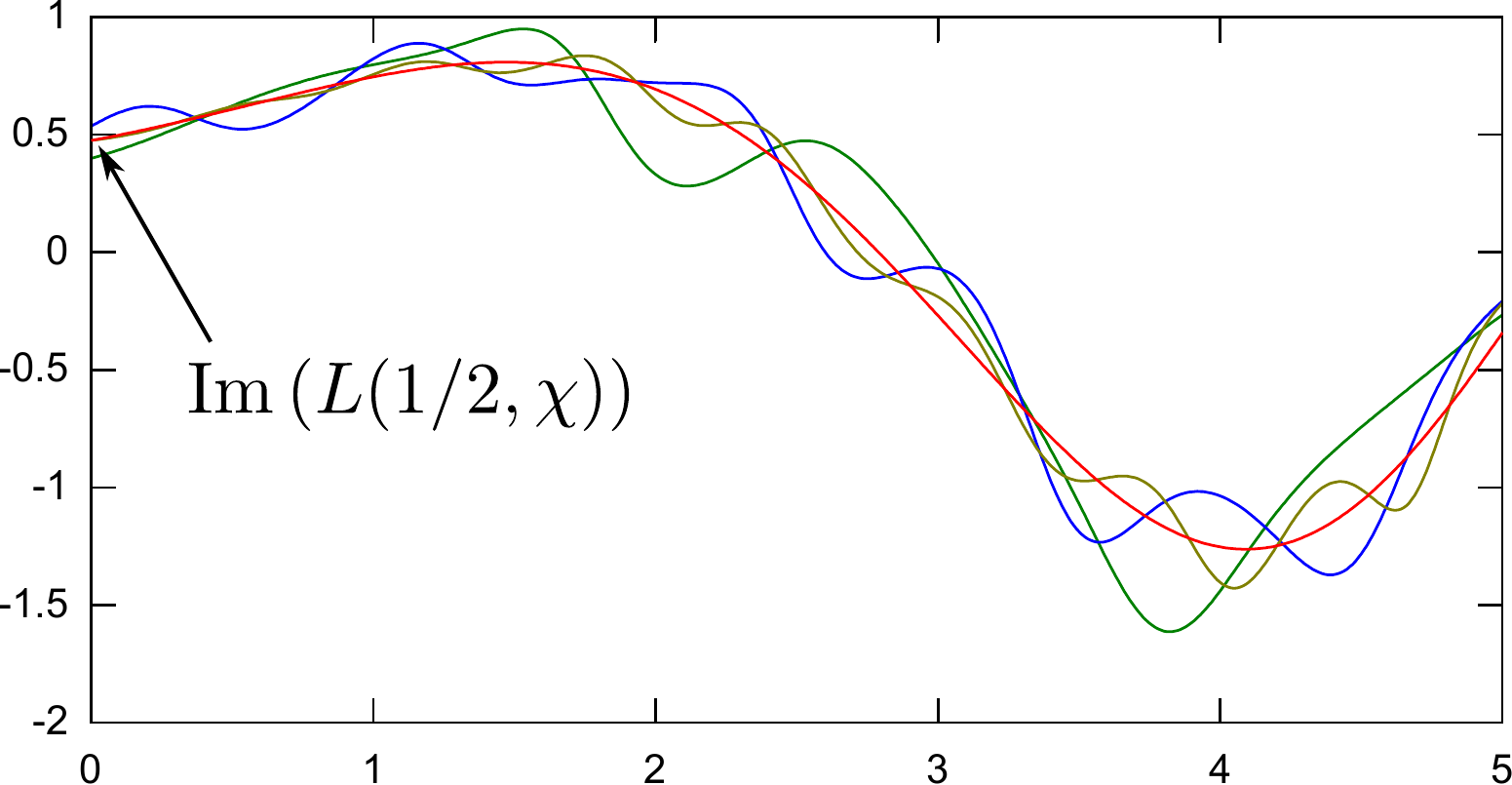} 
 \end{center}
 \caption{Real part (left) and imaginary part (right) of
 $L_x(1/2+it,\chi_{7a})$}
 \label{conv_mod7a}
\end{figure}

\begin{figure}[h]
 \begin{center}
  \includegraphics[height=10em]{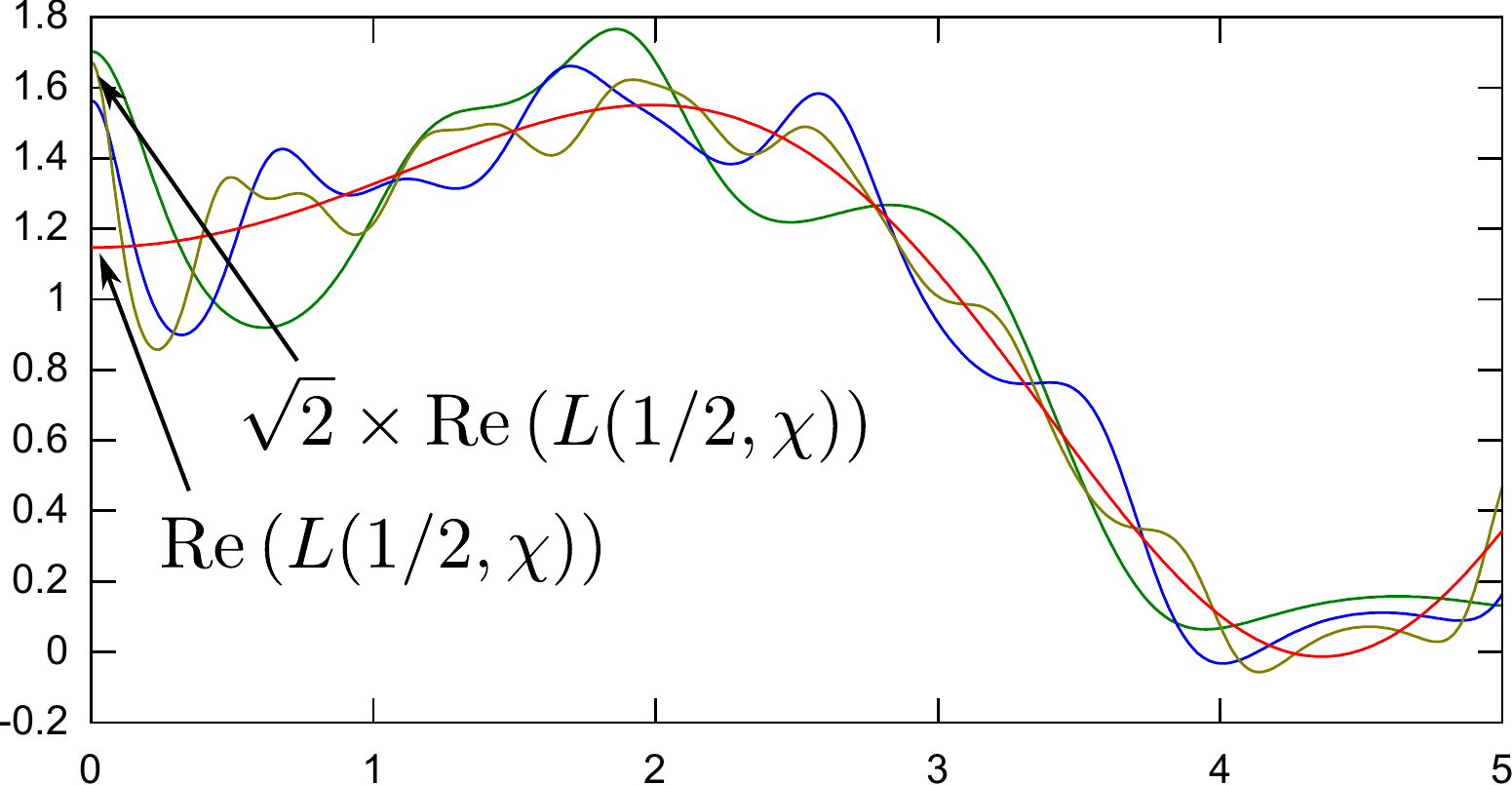} \qquad
  \includegraphics[height=10em]{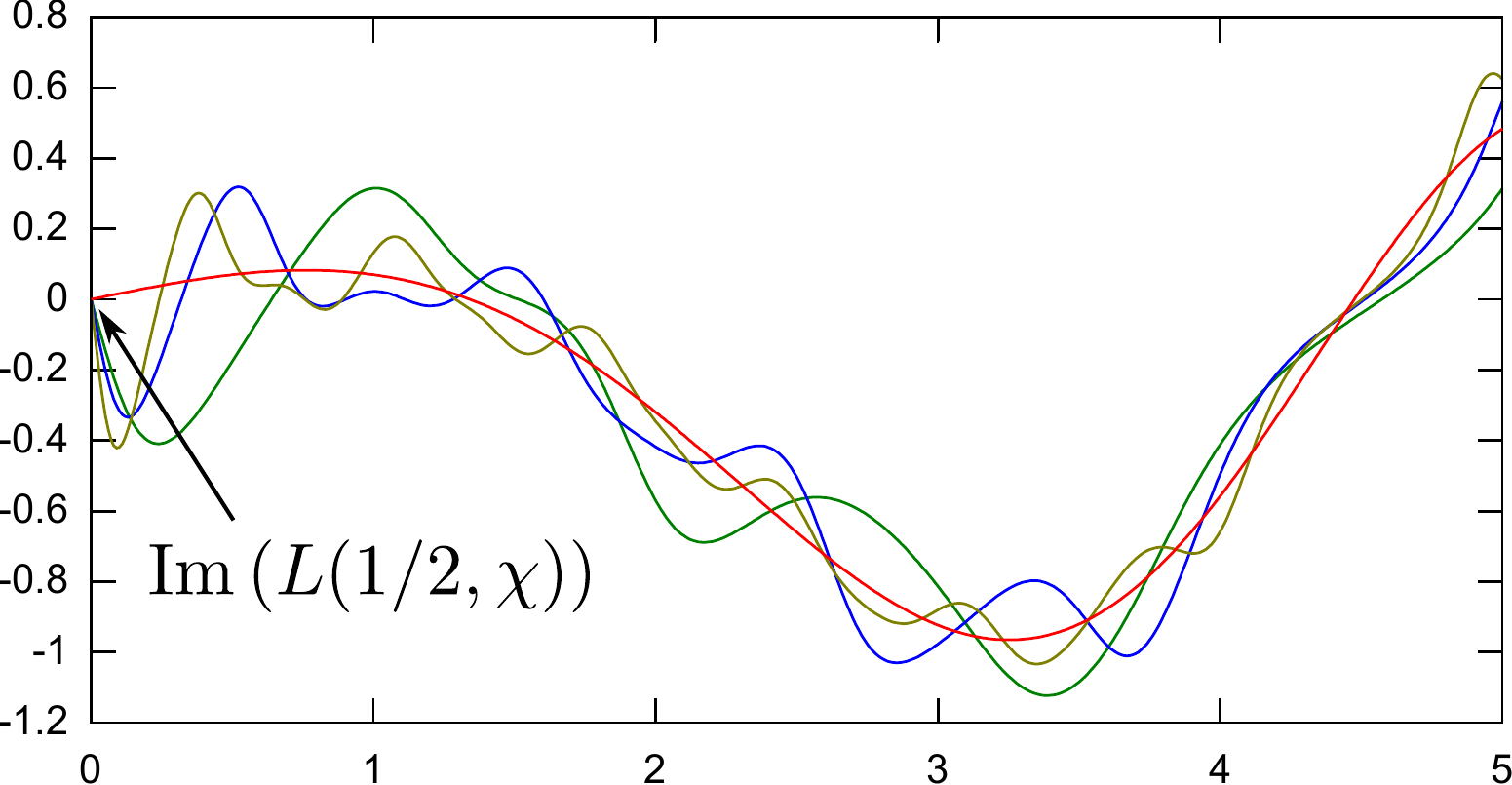} 
 \end{center}
 \caption{Real part (left) and imaginary part (right) of
 $L_x(1/2+it,\chi_{7b})$}
 \label{conv_mod7b}
\end{figure}

\begin{figure}[h]
 \begin{center}
  \includegraphics[height=10em]{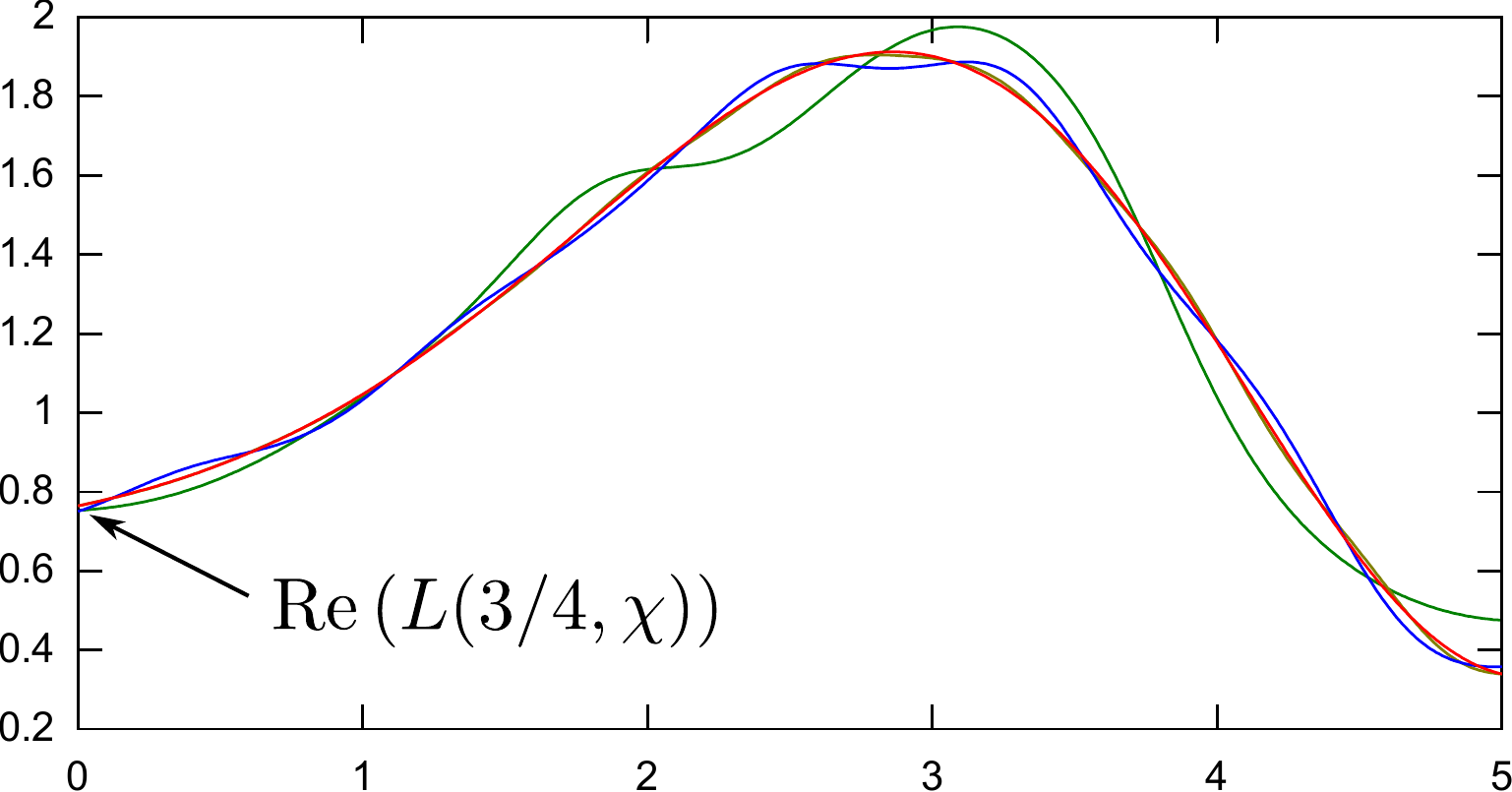} \qquad
  \includegraphics[height=10em]{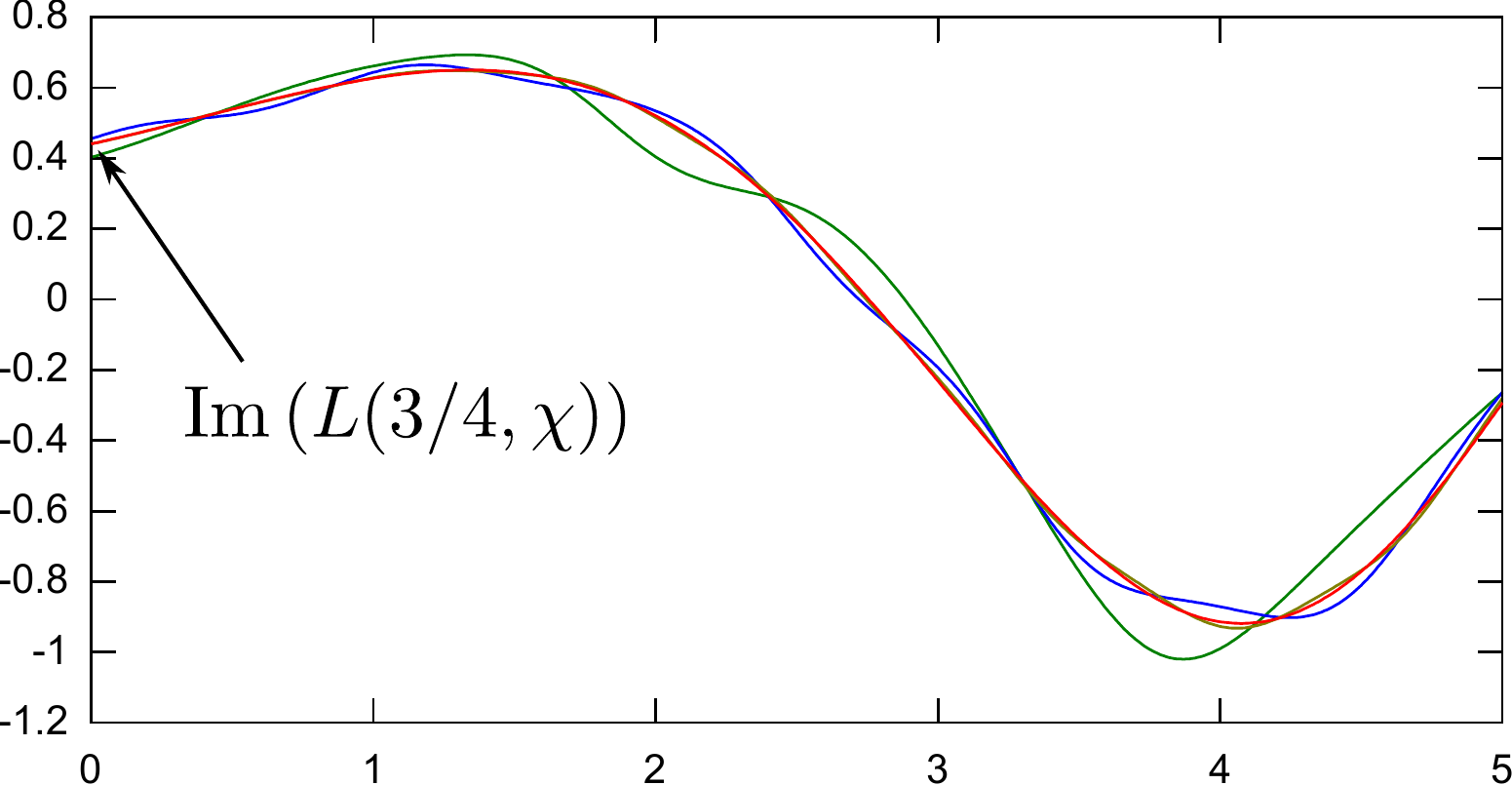} 
 \end{center}
 \caption{Real part (left) and imaginary part (right) of
 $L_x(3/4+it,\chi_{7a})$}
 \label{conv075_mod7a}
\end{figure}

\begin{figure}[h]
 \begin{center}
  \includegraphics[height=10em]{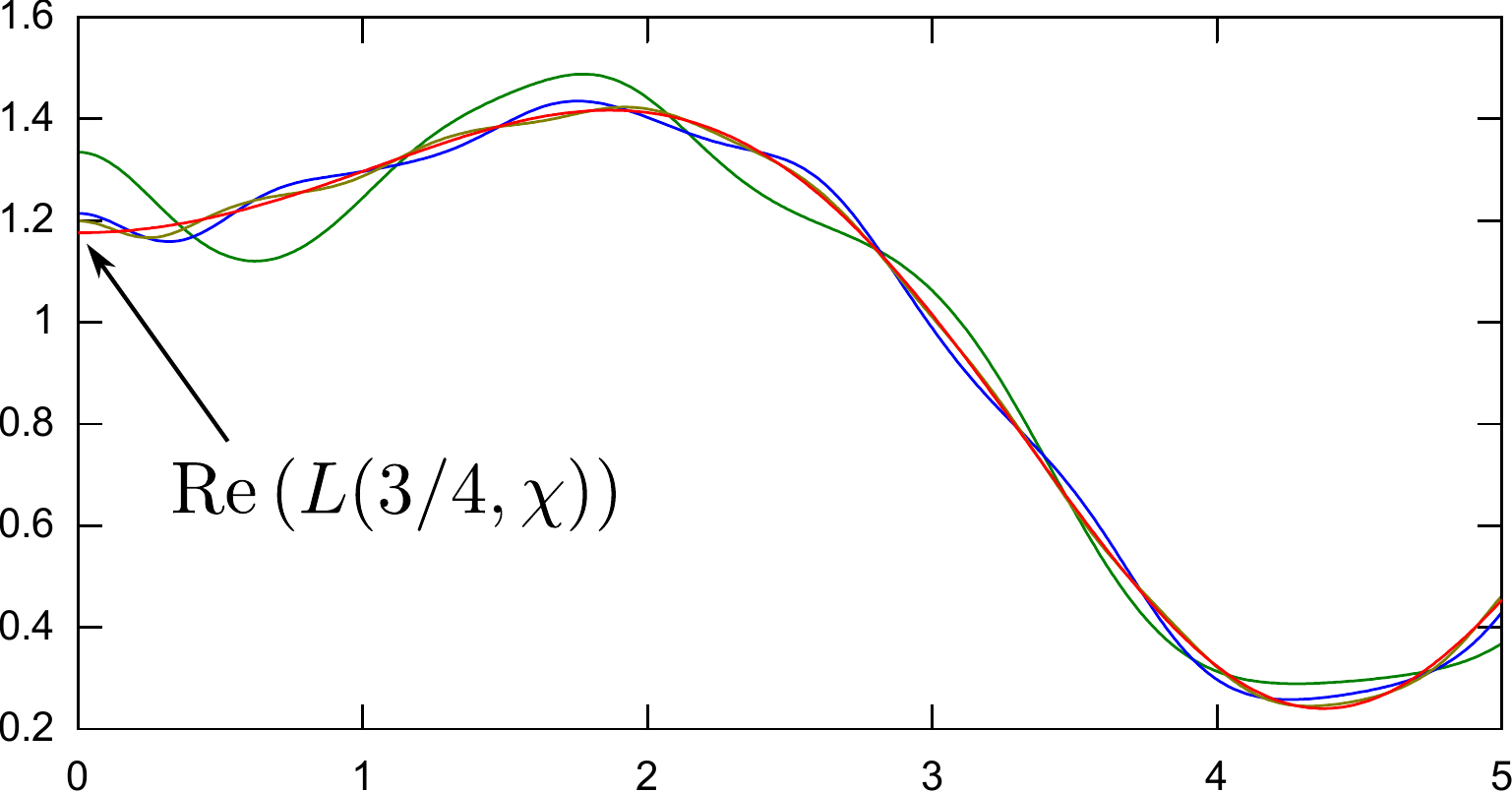} \qquad
  \includegraphics[height=10em]{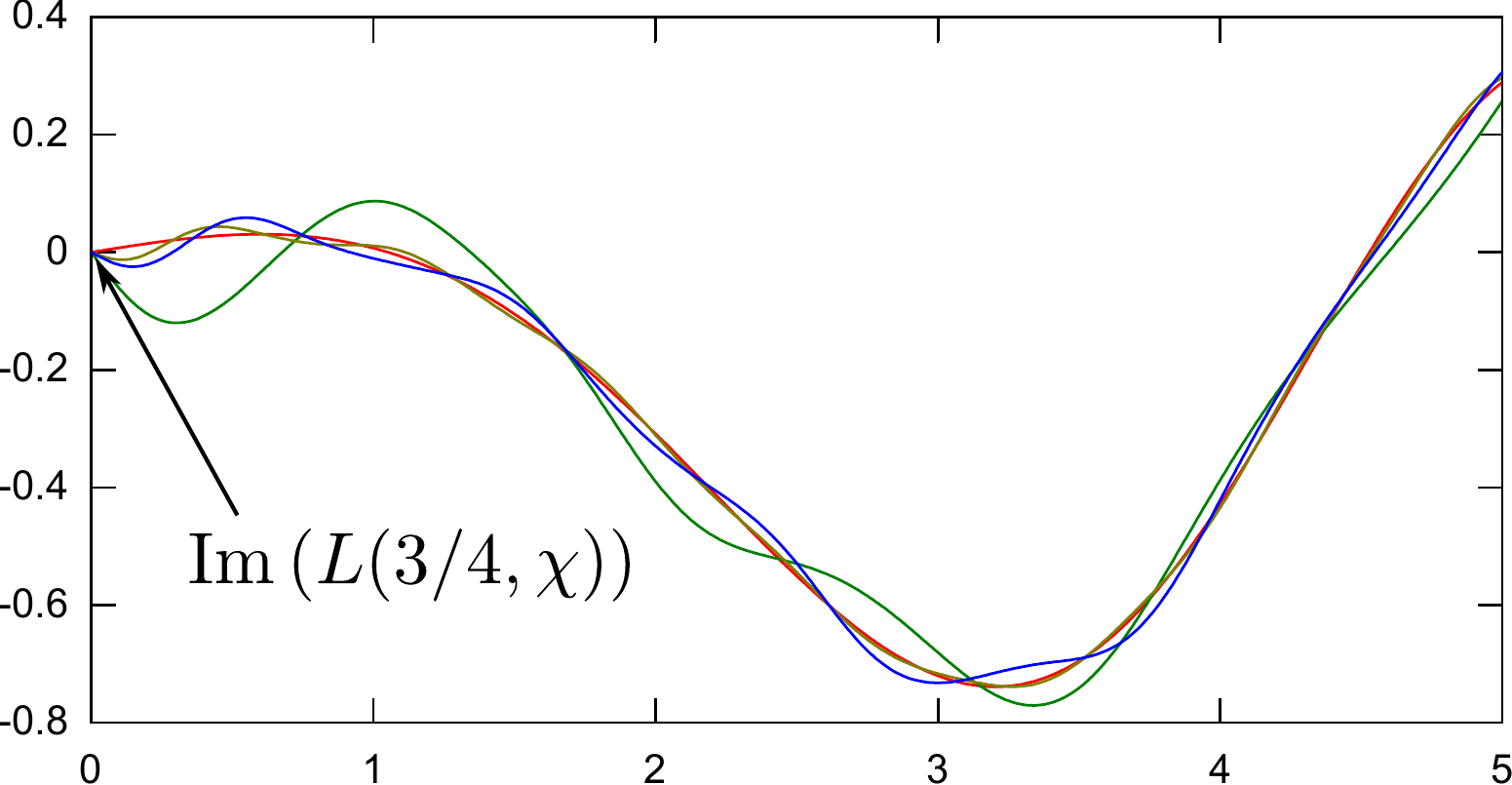} 
 \end{center}
 \caption{Real part (left) and imaginary part (right) of
 $L_x(3/4+it,\chi_{7b})$}
 \label{conv075_mod7b}
\end{figure}

\begin{figure}[h]
 \begin{center}
  \includegraphics[height=10em]{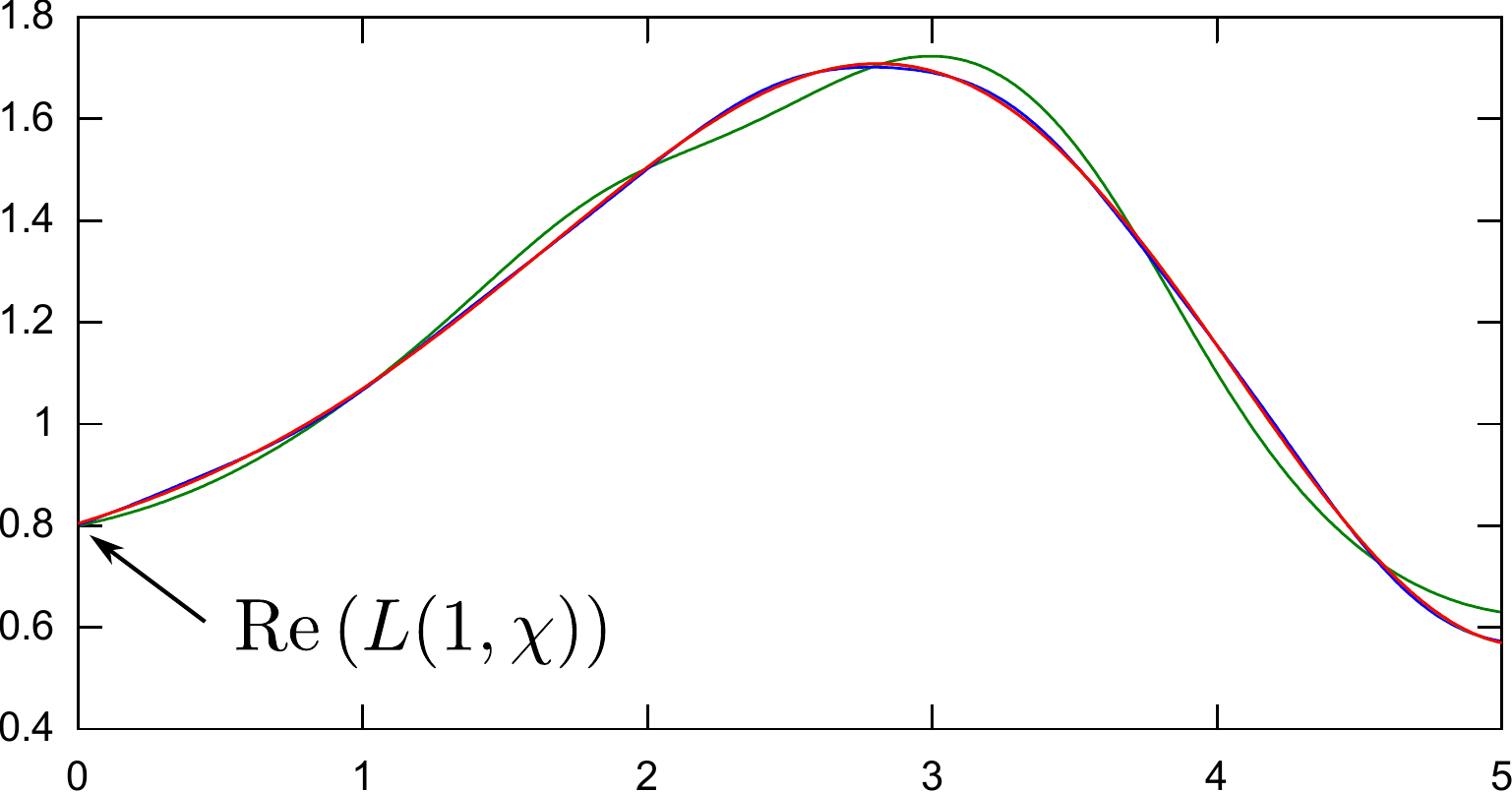} \qquad
  \includegraphics[height=10em]{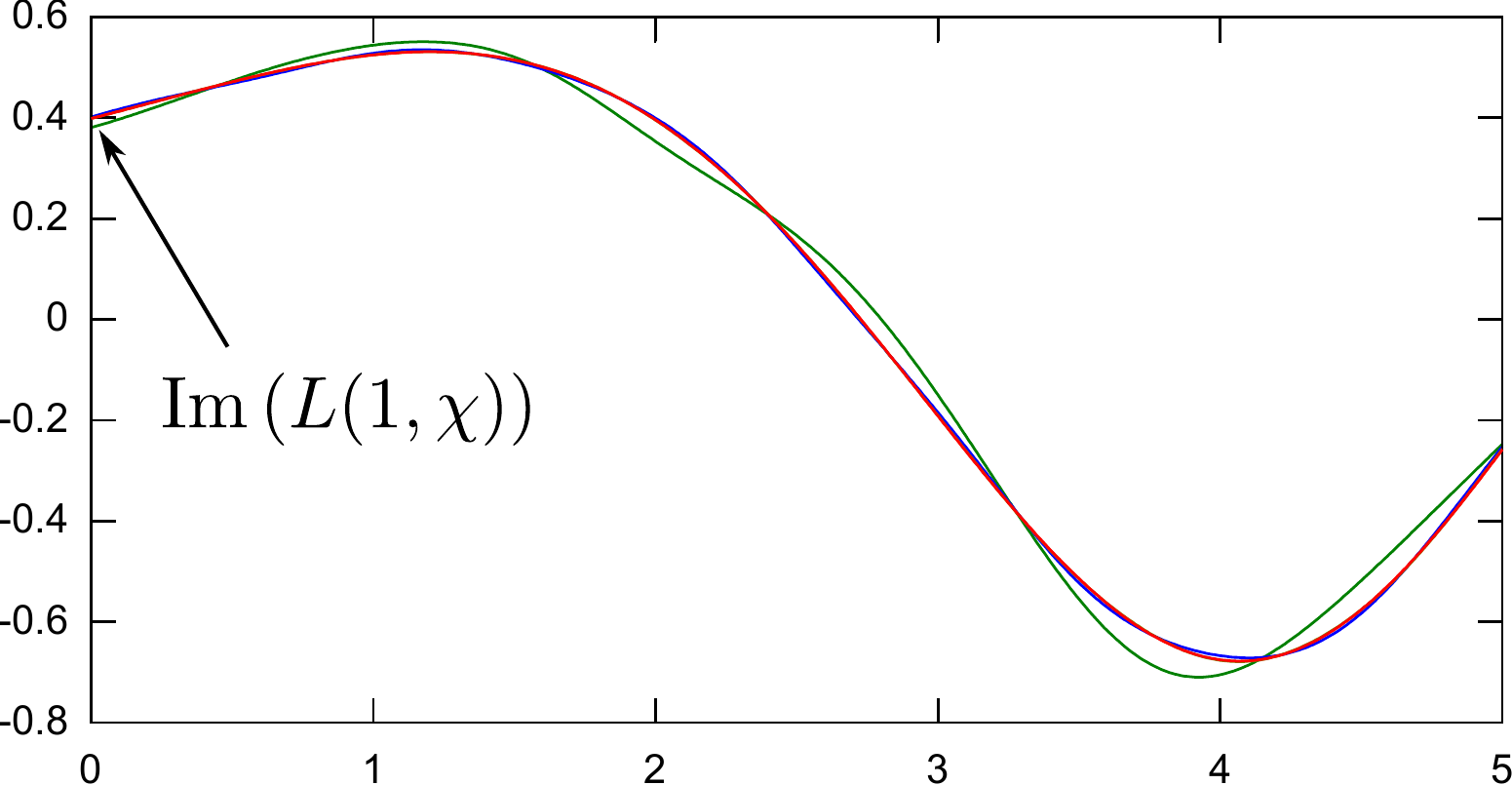} 
 \end{center}
 \caption{Real part (left) and imaginary part (right) of
 $L_x(1+it,\chi_{7a})$}
 \label{conv100_mod7a}
\end{figure}

\begin{figure}[h]
 \begin{center}
  \includegraphics[height=10em]{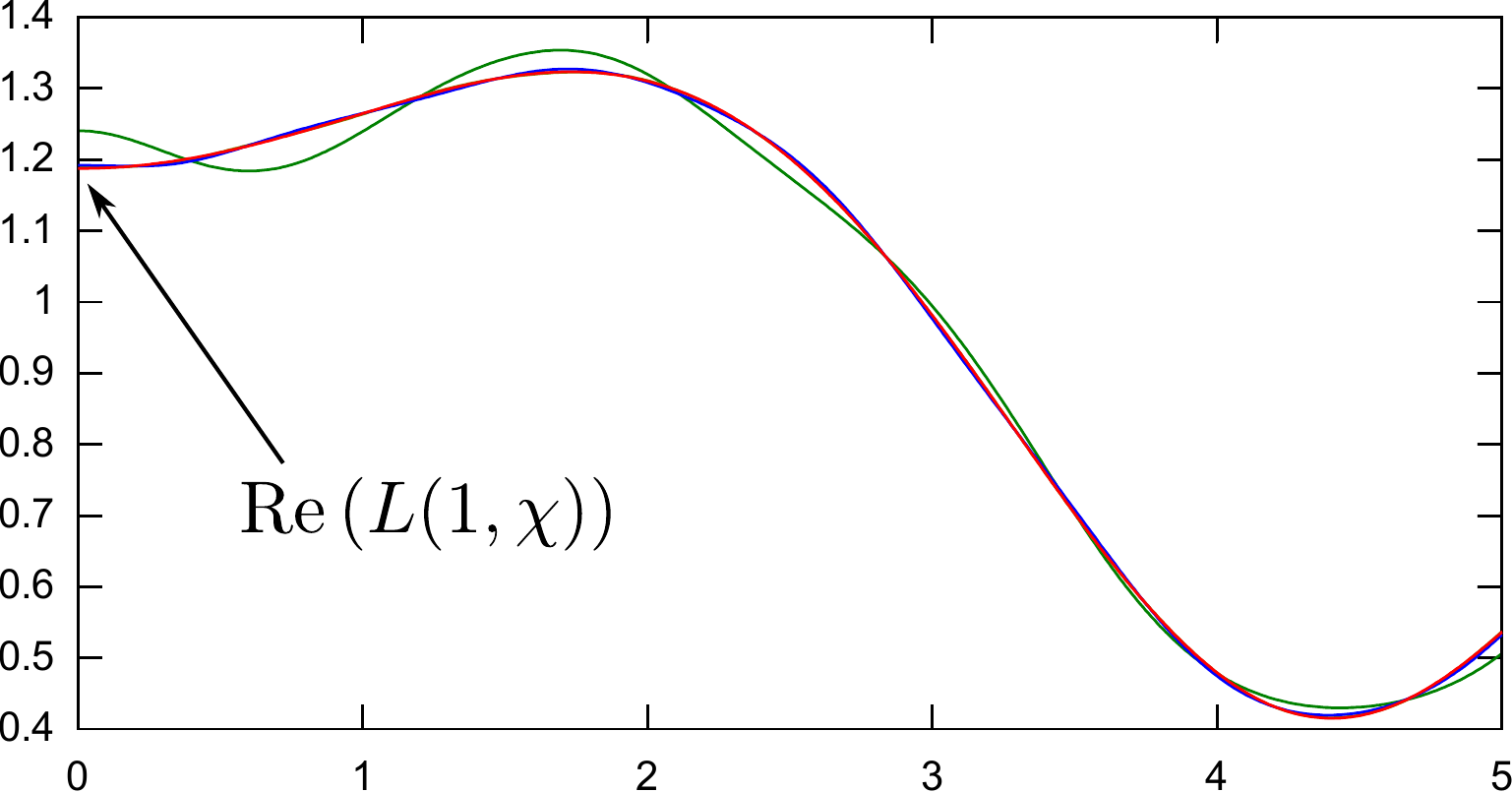} \qquad
  \includegraphics[height=10em]{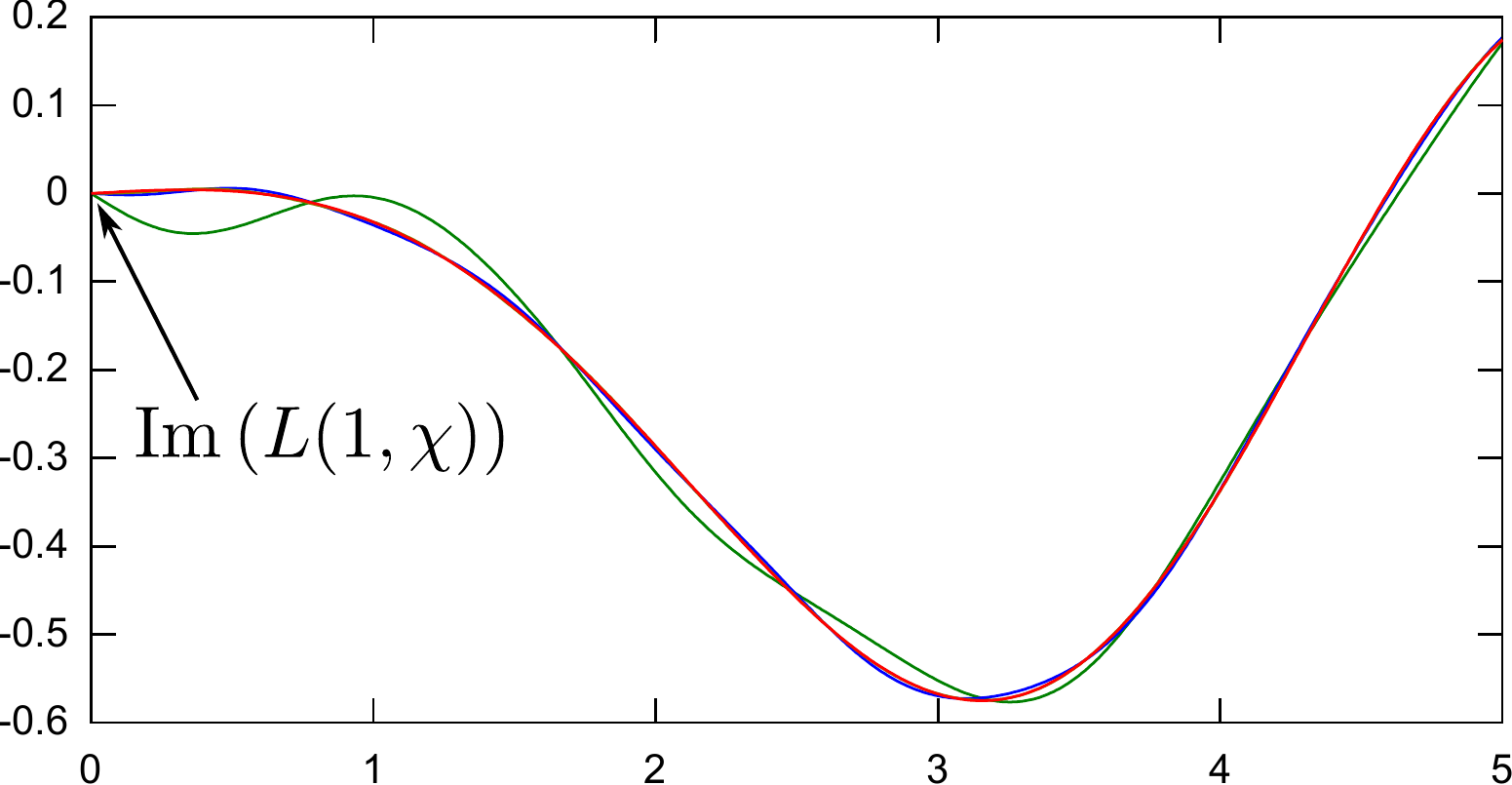} 
 \end{center}
 \caption{Real part (left) and imaginary part (right) of
 $L_x(1+it,\chi_{7b})$}
 \label{conv100_mod7b}
\end{figure}

Denote by $p_n$ the $n$-th prime number.
Figures~\ref{conv_mod7a}, \ref{conv_mod7b}, \ref{conv075_mod7a},
\ref{conv075_mod7b}, \ref{conv100_mod7a} and \ref{conv100_mod7b}
show the datum for the values 
$$L_{x}\l(\frac12+it,\,\chi\r), \qquad
L_{x}\l(\frac34+it,\,\chi\r), \qquad
L_{x}\l(1+it,\,\chi\r)
$$
for $x=p_{10}$ (green), $x=p_{100}$ (blue), $x=p_{1000}$ (yellow) and $\infty$ (red).
Figures~\ref{conv_mod7a}, \ref{conv075_mod7a} and \ref{conv100_mod7a}
are for $\chi_{7a}$,
and Figures~\ref{conv_mod7b}, \ref{conv075_mod7b} and \ref{conv100_mod7b}
for $\chi_{7b}$.
\remnum{
As $t\to0$, we apparently see that
$L_{x}\l(1/2+it,\,\chi\r)\to L(1/2,\chi)$ for $\chi^2\ne\1$, that
$L_{x}\l(1/2+it,\,\chi\r)\to \sqrt2L(1/2,\chi)$ for $\chi^2=\1$, 
and that $L_{x}\l(3/4+it,\,\chi\r) \to L(3/4,\chi)$,
$L_{x}\l(1+it,\,\chi\r) \to L(1,\chi)$ 
for both cases $\chi^2=\1$ and $\chi^2\not=\1$.}{b}
This supports the DRH (Conjecture \ref{conj2}).

\begin{figure}[h]
 \begin{center}
  \includegraphics[width=19em]{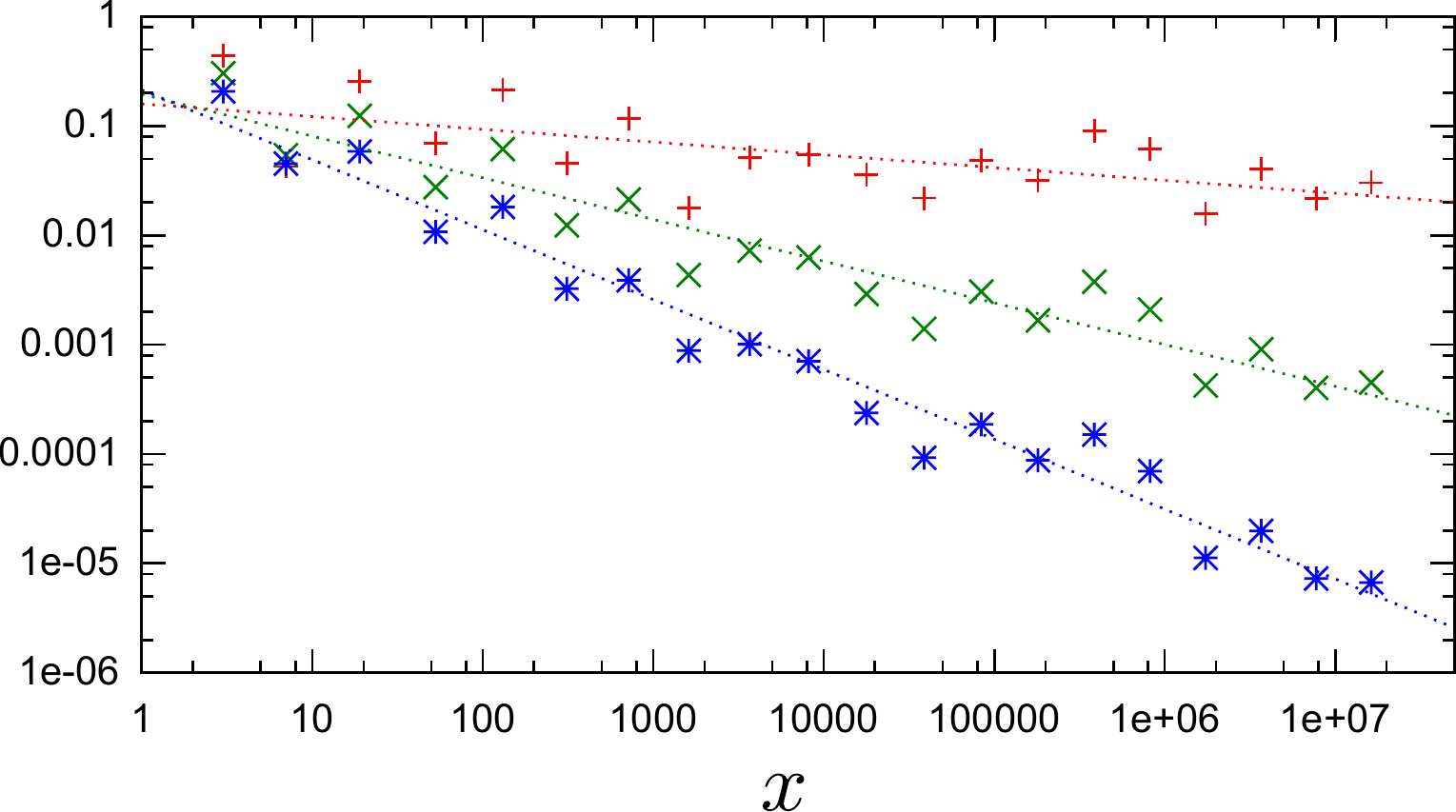} \quad
  \includegraphics[width=19em]{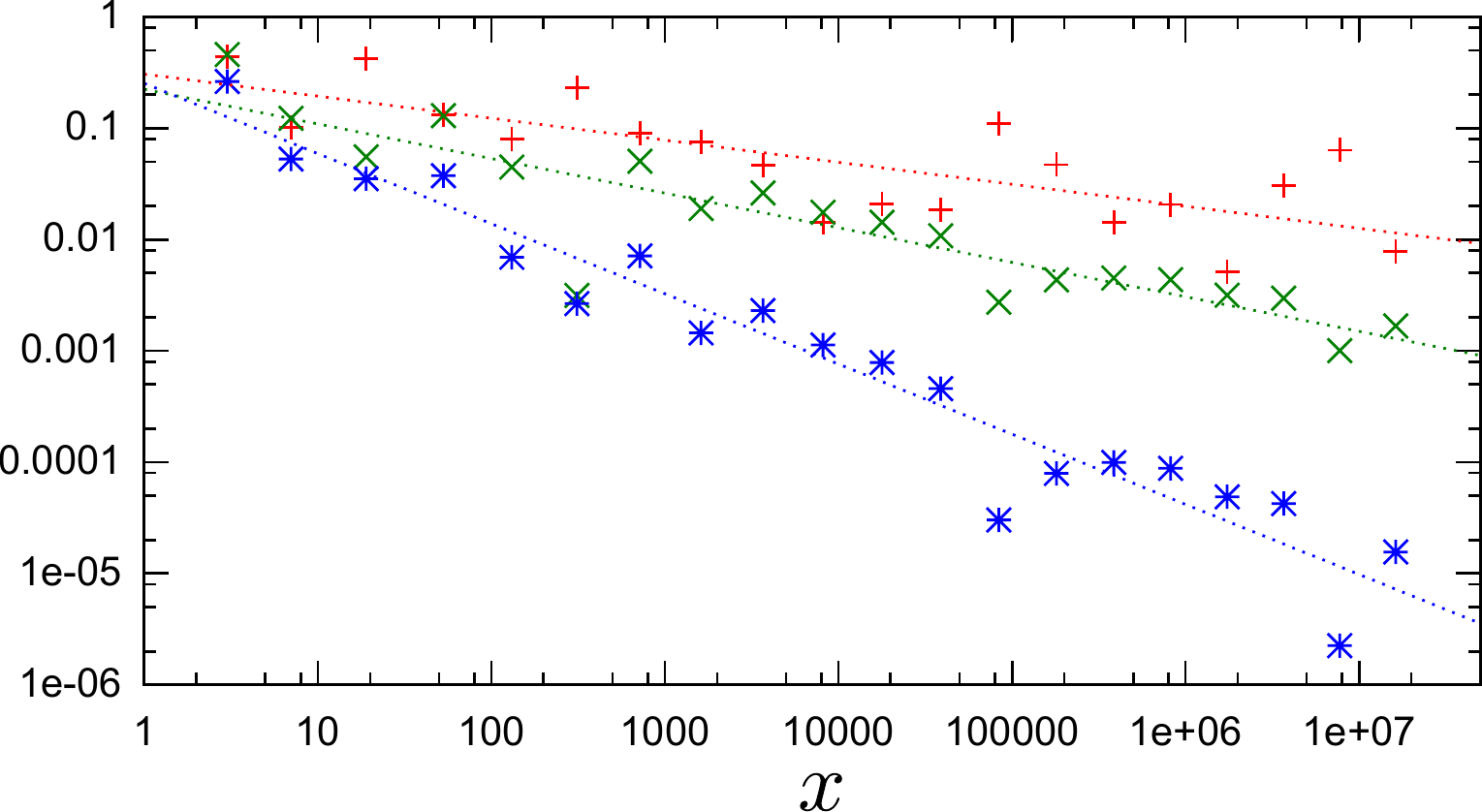} \quad
 \end{center}
 \caption{$\delta L_x(s,\chi)$ for $s=1/2$ (red),
 $s=3/4$ (green) and $s=1$ (blue) with $\chi_{7a}$ (left)
 and $\chi_{7b}$ (right)}
 \label{error_func}
\end{figure}

We introduce the following error function in order to estimate the speed
of convergence for $L_x(s,\chi)$:
\[
  \delta L_x(s,\chi) =
  \begin{cases}
  \left|
   \frac{L_x(s,\chi) - \sqrt{2} L(s,\chi)}{\sqrt{2} L(s,\chi)}
  \right| & (s = 1/2 \ and \ \chi^2=\1)  \\[8pt]
  \left|
   \frac{L_x(s,\chi) - L(s,\chi)}{L(s,\chi)}
  \right| & (otherwise)  
  \end{cases} \, .
\]
Figure~\ref{error_func} shows the values of $\delta L_x(s,\chi)$.
When we approximate the error function as $\delta L_x(s,\chi) \sim
x^{-\alpha}$, the exponents are determined 
\remnum{so that they fit}{c}
the numerical results (Table~\ref{exponent}).
We see the speed 
\remnum{
of convergence becomes faster as $s$ gets larger, if $s$ is real.}{d}

\begin{table}[t]
 \begin{center}
    \begin{tabular}{ccccc} \hline \hline
    $s$ && $\alpha \ (\chi_{7a})$ && $\alpha \ (\chi_{7b})$ \\ \hline
     $1/2$ && 0.1167 && 0.1978 \\ 
     $3/4$ && 0.3814 && 0.3106 \\ 
     $1$ && 0.6389 && 0.6302 \\ \hline \hline
    \end{tabular}
 \end{center}
 \caption{Exponents of $\delta L_{x}(s,\chi)\sim
 x^{-\alpha}$ for $\chi_{7a}$ and $\chi_{7b}$.}
 \label{exponent}
\end{table}

\section{Finite Size Scaling}

\remnum{
In this section, we show another special feature that $L_x(s,\chi)$ has.
Since $L_x(s,\chi)$ is a finite Euler product, it obviously has no zeros
on the critical line.
Nevertheless, $L_x(s,\chi)$ gives a certain sequence of complex numbers,
which seemingly grows up to the nontrivial zeros of $L(s,\chi)$, as
$x\to\infty$. 
In other words, the finite partial Euler product $L_x(s,\chi)$ already
``knows'' the nontrivial zeros of $L(s,\chi)$.}{46}

\begin{figure}[h]
 \begin{minipage}[c]{0.5\textwidth}
  \begin{center}
   \includegraphics[width=20em]{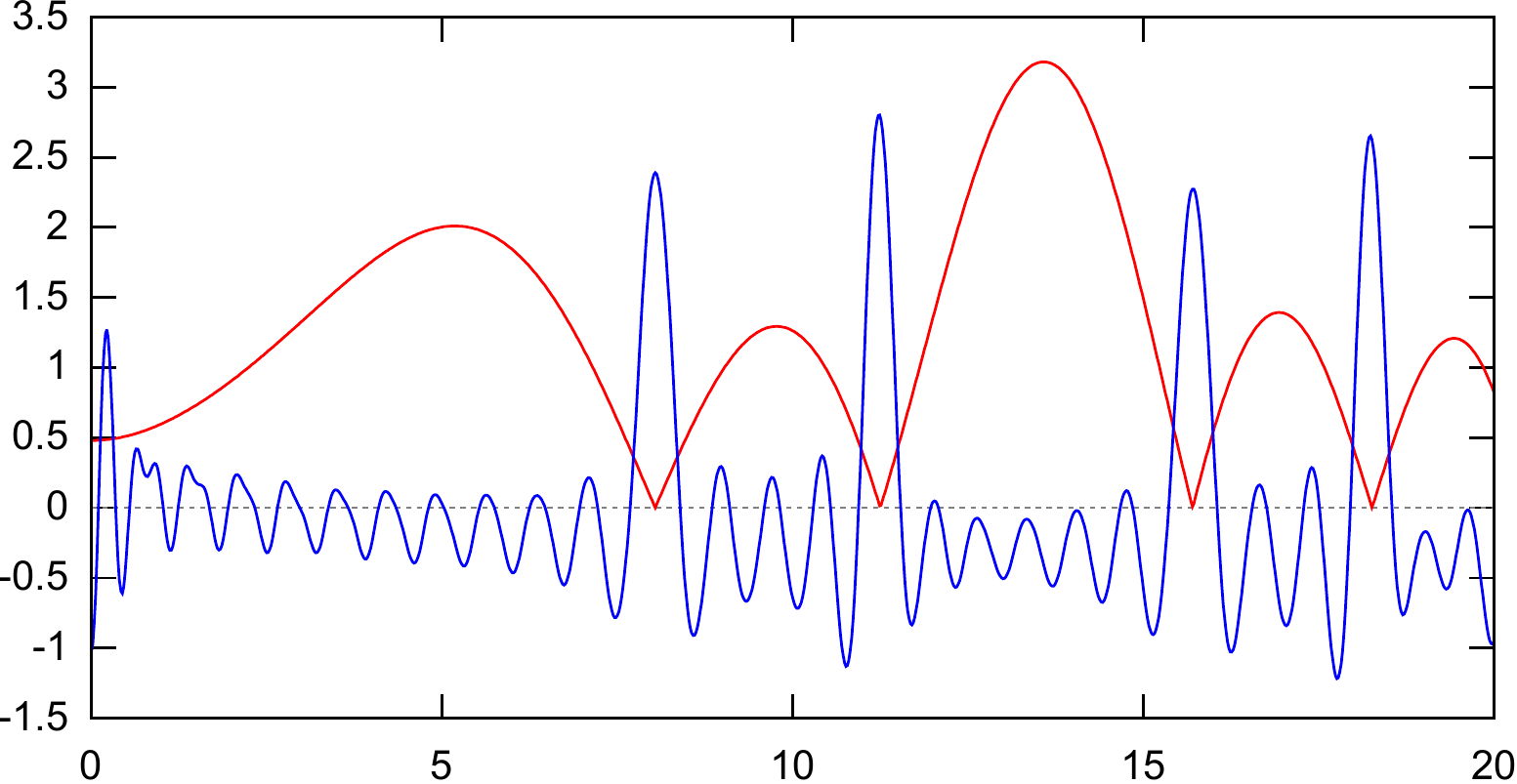}
  \end{center}
  \caption{$\rho_x(t)$ for $\chi_3$}
  \label{Figure3}
 \end{minipage}

 \vspace{1em}

 \begin{minipage}[c]{0.5\textwidth}
  \begin{center}
   \includegraphics[width=20em]{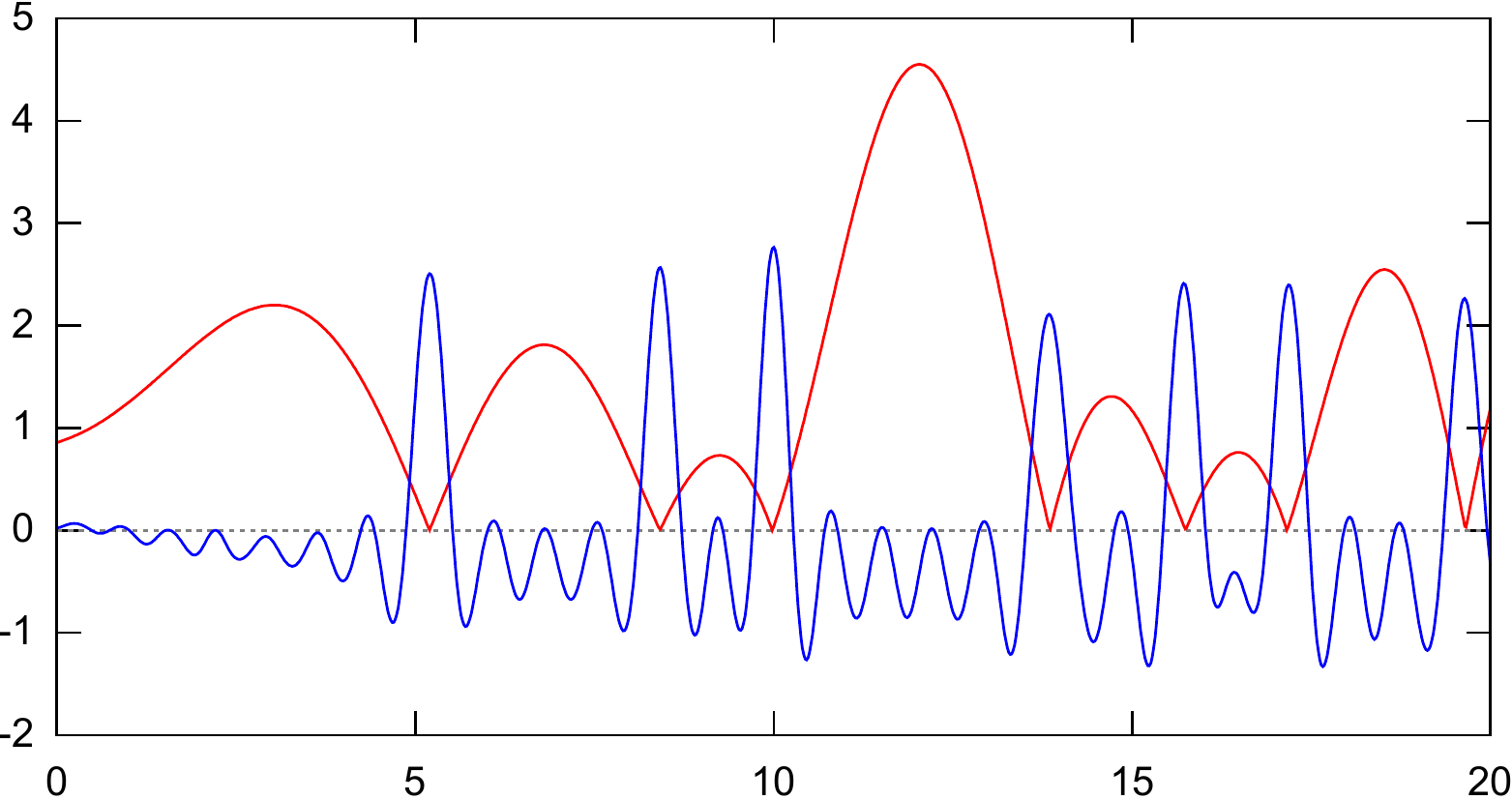}
  \end{center}
  \caption{$\rho_x(t)$ for $\chi_{7a}$}
  \label{Figure4}
 \end{minipage}%
 \begin{minipage}[c]{0.5\textwidth}
  \begin{center}
   \includegraphics[width=20em]{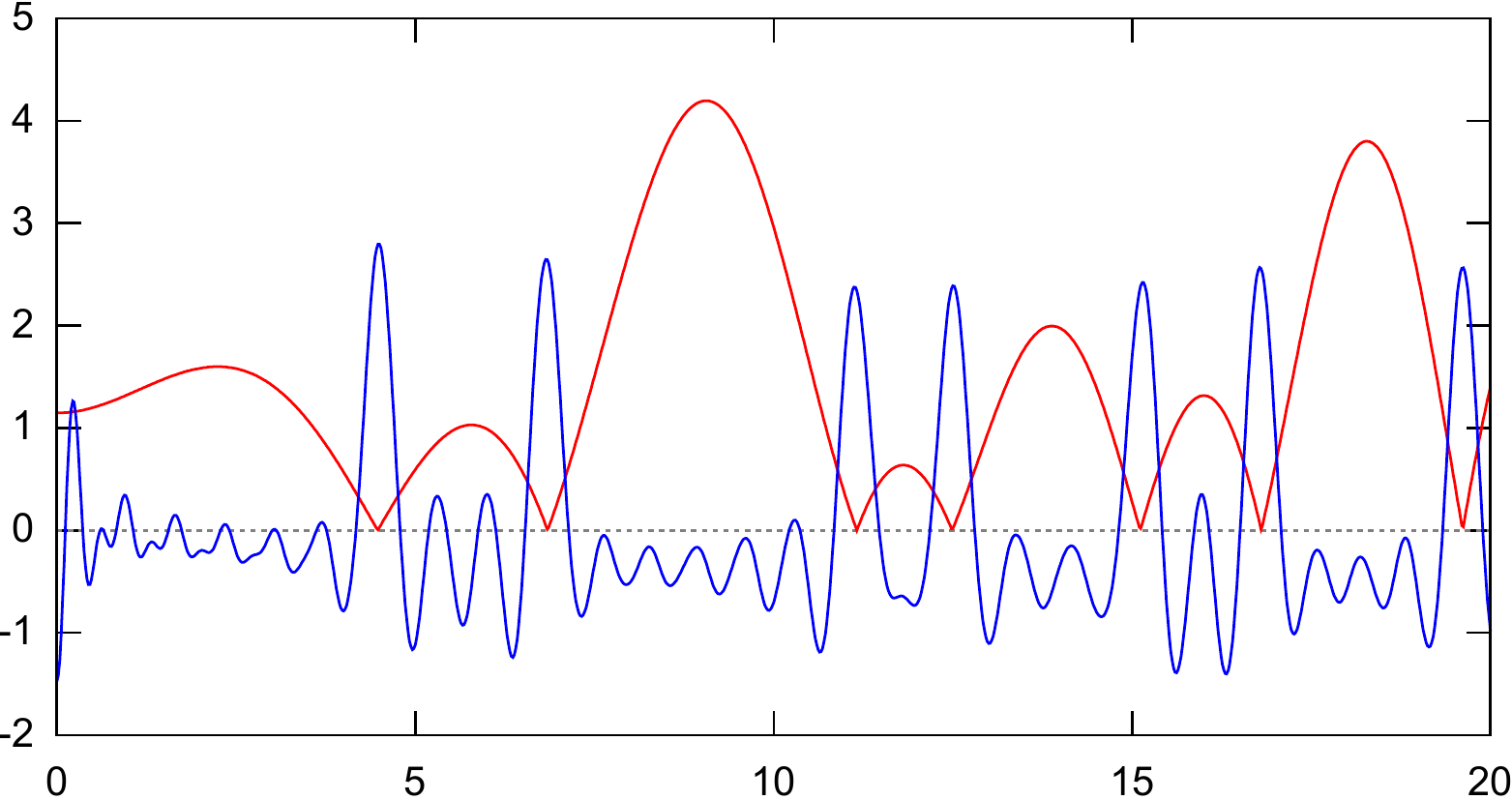}
  \end{center}
  \caption{$\rho_x(t)$ for $\chi_{7b}$}
  \label{Figure5}
 \end{minipage}
\end{figure}

In Figures~\ref{Figure3}, \ref{Figure4} and \ref{Figure5}, the blue
curves show the values
\begin{equation}
 \rho_x(t)=\frac1\pi\Im\frac d{dt}\log L_{x}\l(\frac12+it,\,\chi\r)
  \label{density_func}
\end{equation}
with $x=p_{1000}$
for $\chi_{3}$, $\chi_{7a}$, $\chi_{7b}$, respectively.
The red curves are $|L\l(\frac12+it,\,\chi\r)|$.
This function~(\ref{density_func}) is an analog of the eigenvalue density
function in random matrix theory.
\remnum{The Riemann zeta}{e} function on the critical line $s=1/2+it$
can be seen as a characteristic polynomial of a certain infinite
dimensional matrix~\cite{KS,BH}:
\remnum{
With the Riemann-Siegel theta function
\[
  \vartheta(t) = \Im \log \Gamma\left(\frac{it}{2}+\frac{1}{4}\right)
  - \frac{t}{2} \log \pi \, ,
\]
the function $Z(t)=e^{i\vartheta(t)}\zeta\left(\frac{1}{2}+it\right)$ turns
out to be real.
This is because the completed $\zeta$-function 
\[
 \xi(s) = \frac{s(s-1)}{2} \pi^{-s/2} \Gamma \left(\frac{s}{2}\right) \zeta(s)
\]
is real on $\Re(s)=1/2$ due to the functional equation
$\xi(s)=\xi(1-s)$. 
Dirichlet $L$-functions also have similar representations.
The real function $Z(t)$ changes its signature at nontrivial zeros of the
Riemann zeta function. Thus $Z(t)$ is expressed as
a regularized product
\[
   \prod_{j=1}^\infty \raise1.5ex\hbox{\scriptsize reg} \, ( t - t_j ) 
\]
where $t_j$ satisfies $\zeta\left(\frac12+it_j\right)=0$.
This means the argument of $Z(t)$ jumps by $\pi$ at the zeros.}{38}
Therefore when we define the density function of the nontrivial zeros on the
critical line as
\begin{eqnarray}
 \rho(t) & = & \sum_{j=1}^\infty \delta(t - t_j)
  \nonumber \\
 & = &
  \frac{1}{\pi} \Im
  \sum_{j=1}^\infty \frac{1}{t-t_j}
  \nonumber \\
 & = &
  \frac{1}{\pi} \Im \frac{d}{dt} \log \prod_{j=1}^\infty
  \raise1.5ex\hbox{\scriptsize reg} \, (t - t_j) ,
  \nonumber
\end{eqnarray}
the function (\ref{density_func}) should converge to this density
function in the limit of $x\to\infty$,
up to the factor coming from $\vartheta(t)$.
\remnum{
Here we simply write the delta function as $\delta(x) = \frac{1}{\pi}
\Im \frac{1}{x}$, which is originally represented as
$\delta(x)=\lim_{\epsilon\to 0^+}\pm\frac{1}{\pi}\Im\frac{1}{x\mp i\epsilon}$.}{38}

Apparently the location of the zeros of $|L(\frac12+it,\,\chi)|$ agrees
to that of the peaks of $\rho_x(t)$ in Figures~\ref{Figure3},
\ref{Figure4} and \ref{Figure5}.
This suggests that 
\remnum{a finite set of first few primes already ``knows''}{f}
the nontrivial zeros of $L(s,\,\chi)$, and that the Euler product 
\remnum{would be}{g}
meaningful beyond the boundary.
We also observe that the blue curve oscillates near $t=0$ if and only if
$\chi^2=\1$.

\begin{figure}[h]
 \begin{center}
  \includegraphics[height=18em]{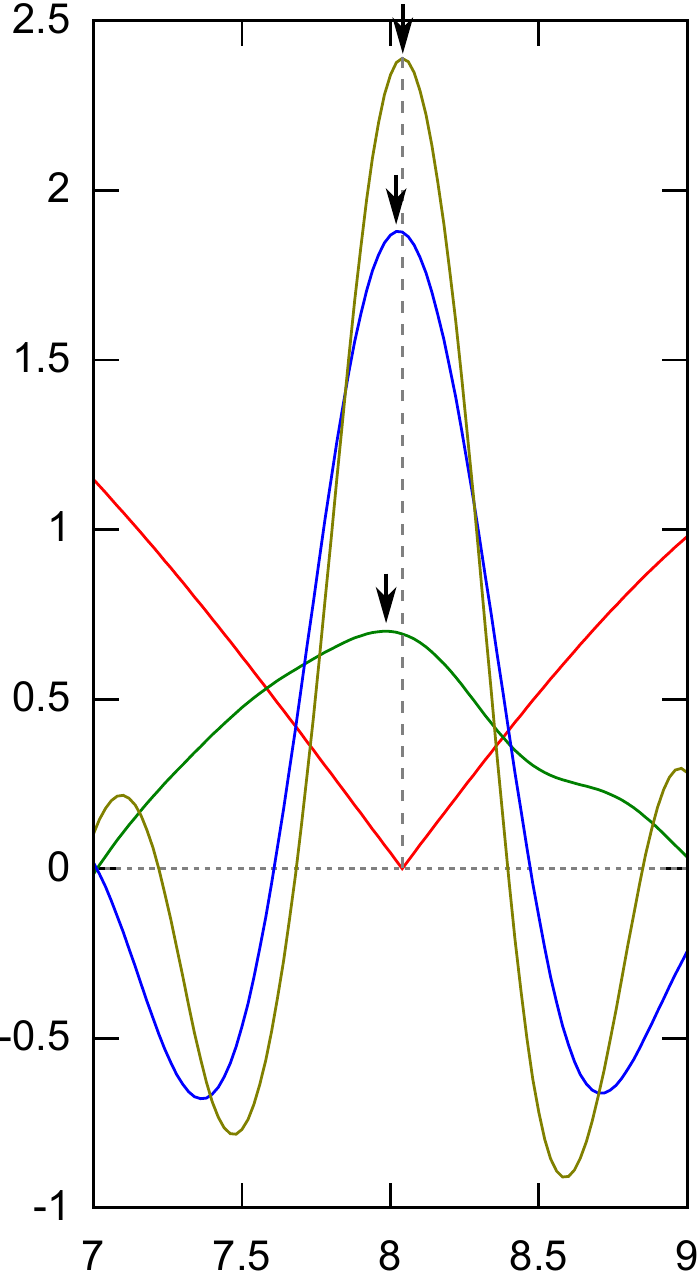} \qquad
  \includegraphics[height=18em]{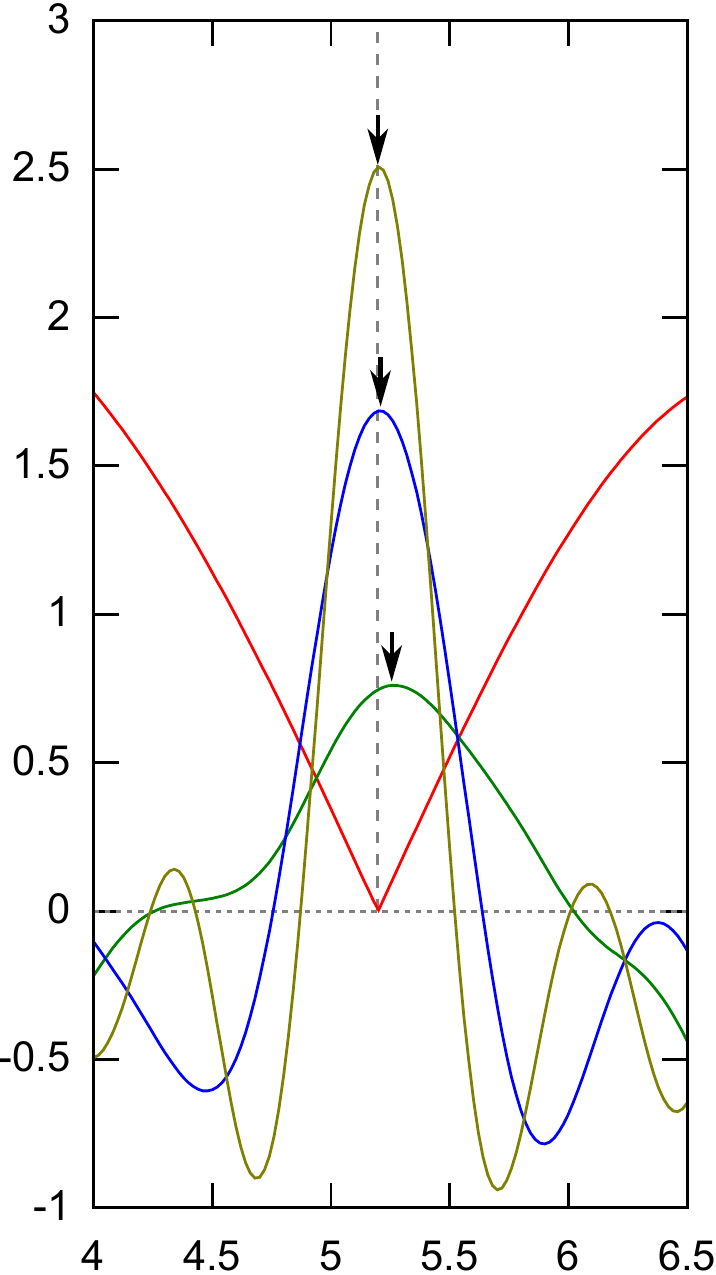} \qquad
  \includegraphics[height=18em]{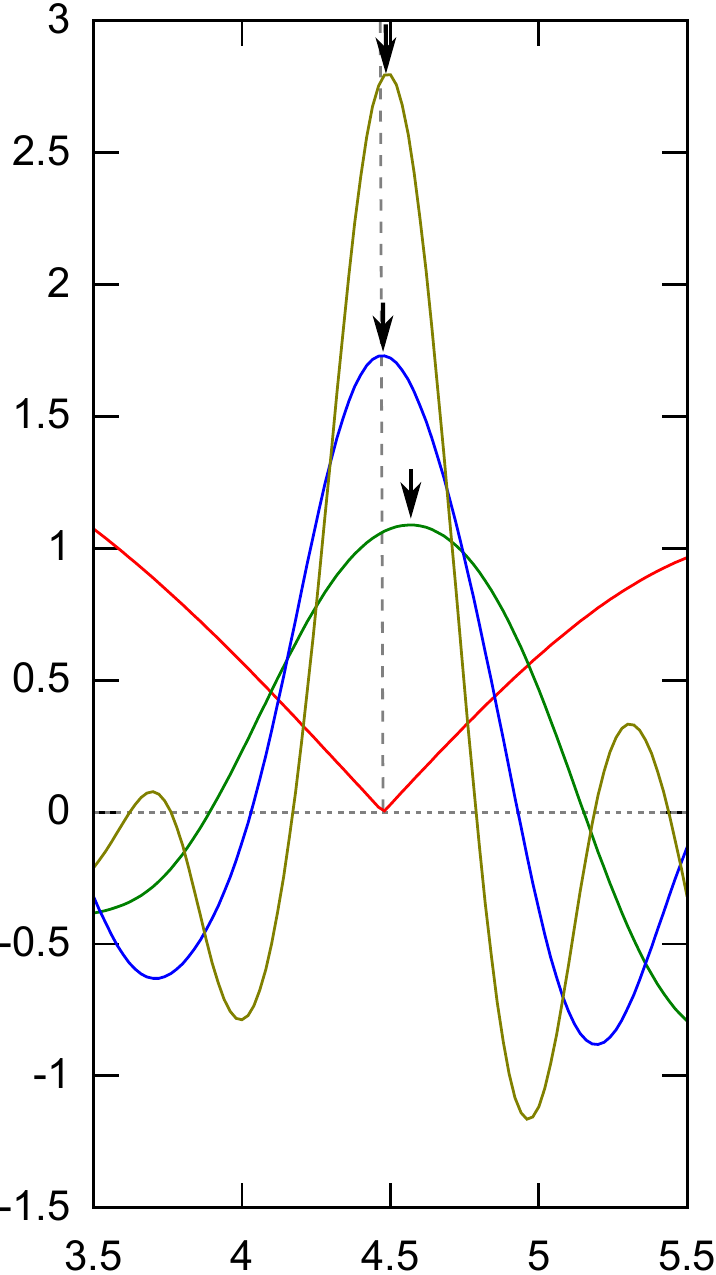} 
 \end{center}
 \caption{Peaks in $\rho(t)$ with the smallest zero for $\chi_3$ (left),
 $\chi_{7a}$ (center) and $\chi_{7b}$ (right)}
 \label{peaks}
\end{figure}

Figure~\ref{peaks} shows how the peaks of $\rho(t)$ with the smallest zero
in Figures~\ref{Figure3}, \ref{Figure4} and \ref{Figure5} get closer to
the zeros of $L(s,\,\chi)$ for $x=p_{10}$ (green), $x=p_{100}$ (blue),
$x=p_{1000}$ (yellow).
We see these peaks getting higher and narrower, and approaching the
Dirac delta function.
This kind of scaling behavior is often found in critical phenomena
associated with some phase transitions.
\remnum{
Especially, in this case, the situation is similar to percolation
theory~\cite{StaufferAharony}.}{40}

\begin{figure}[h]
 \begin{minipage}[c]{0.5\textwidth}
  \begin{center}
     \includegraphics[width=20em]{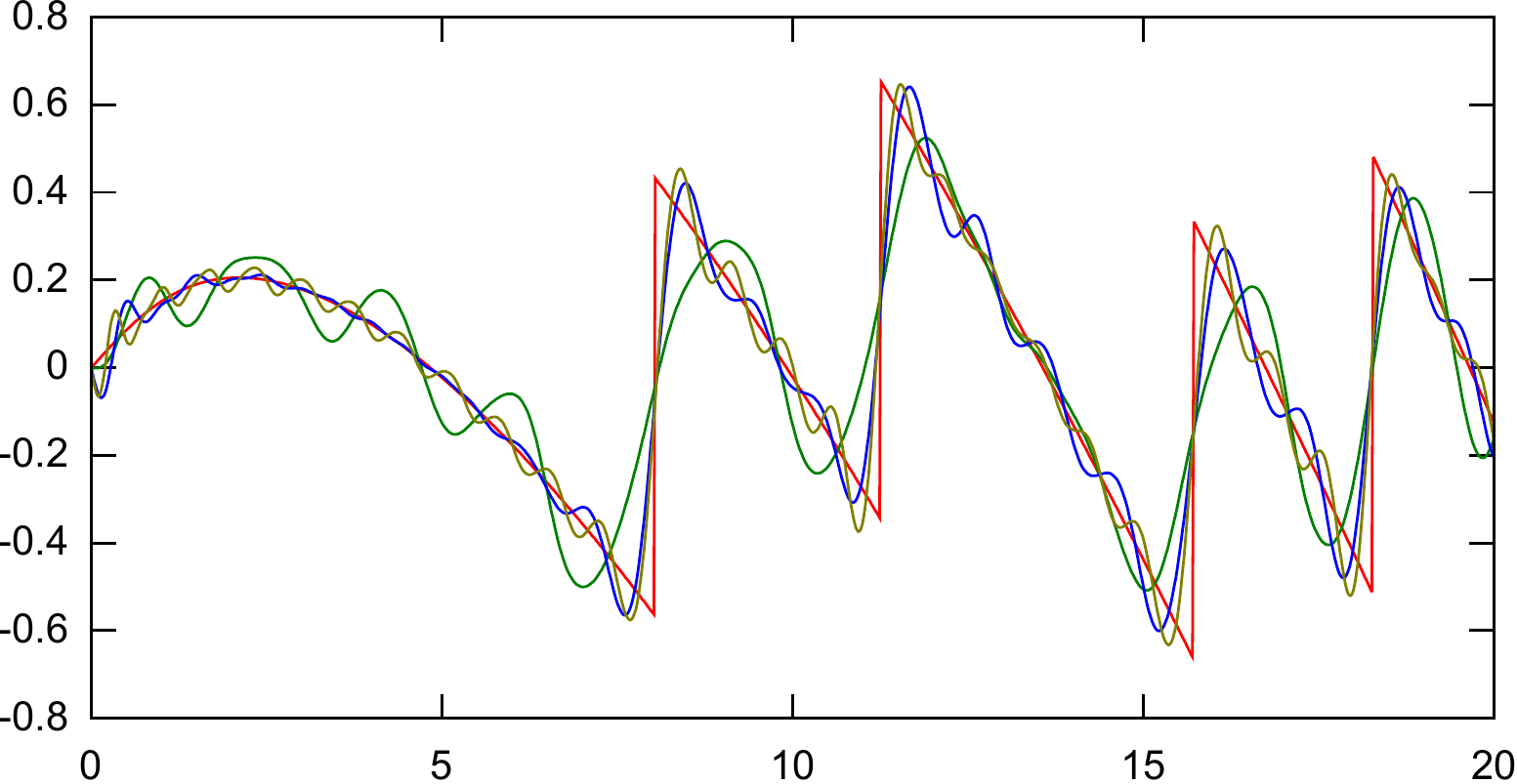}
  \end{center}
  \caption{$R_x(t)$ for $\chi_3$}
  \label{Figure6}
 \end{minipage}

 \vspace{1em}

 \begin{minipage}[c]{0.5\textwidth}
  \begin{center}
     \includegraphics[width=20em]{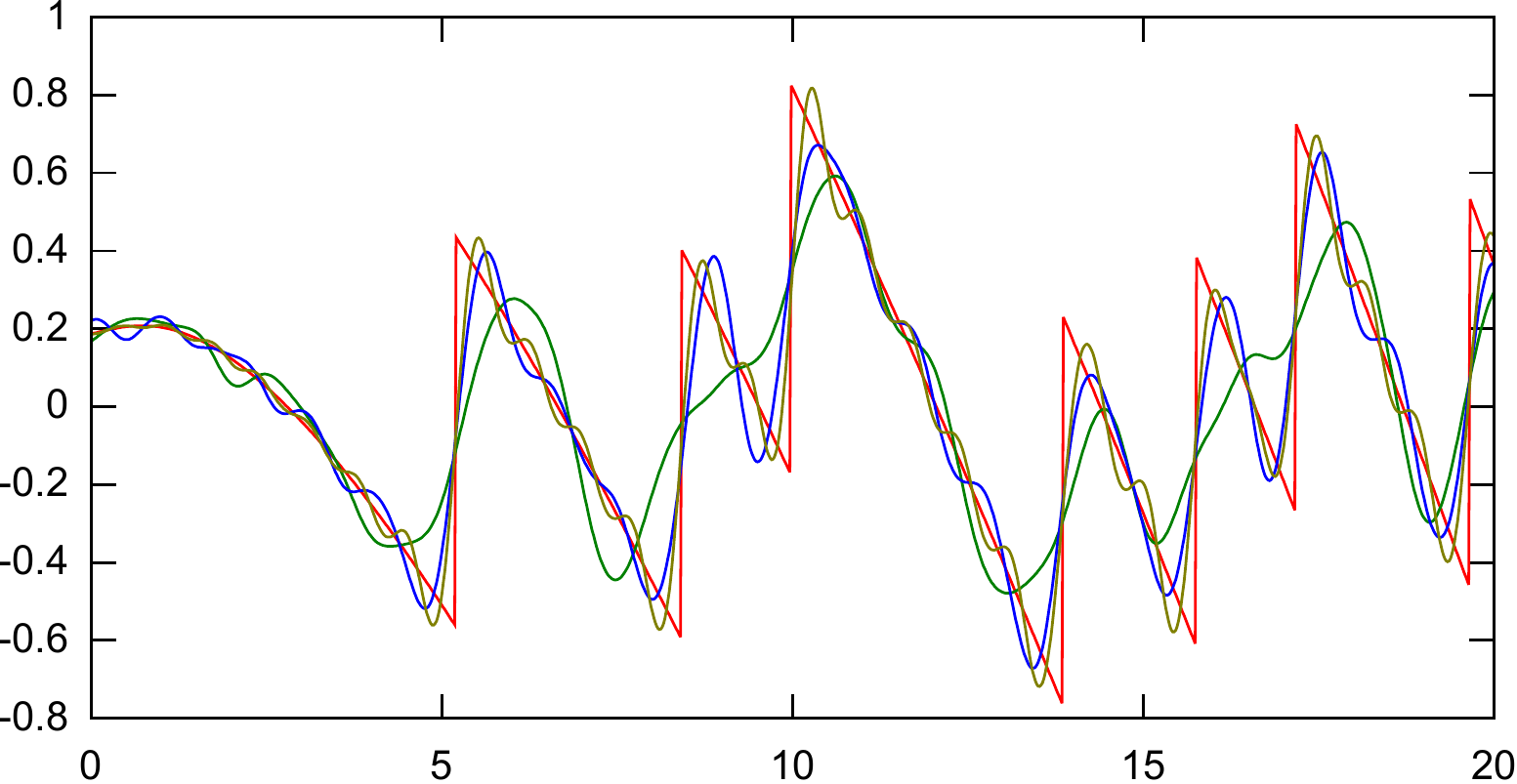}
  \end{center}
  \caption{$R_x(t)$ for $\chi_{7a}$}
  \label{Figure7}
 \end{minipage}%
 \begin{minipage}[c]{0.5\textwidth}
  \begin{center}
     \includegraphics[width=20em]{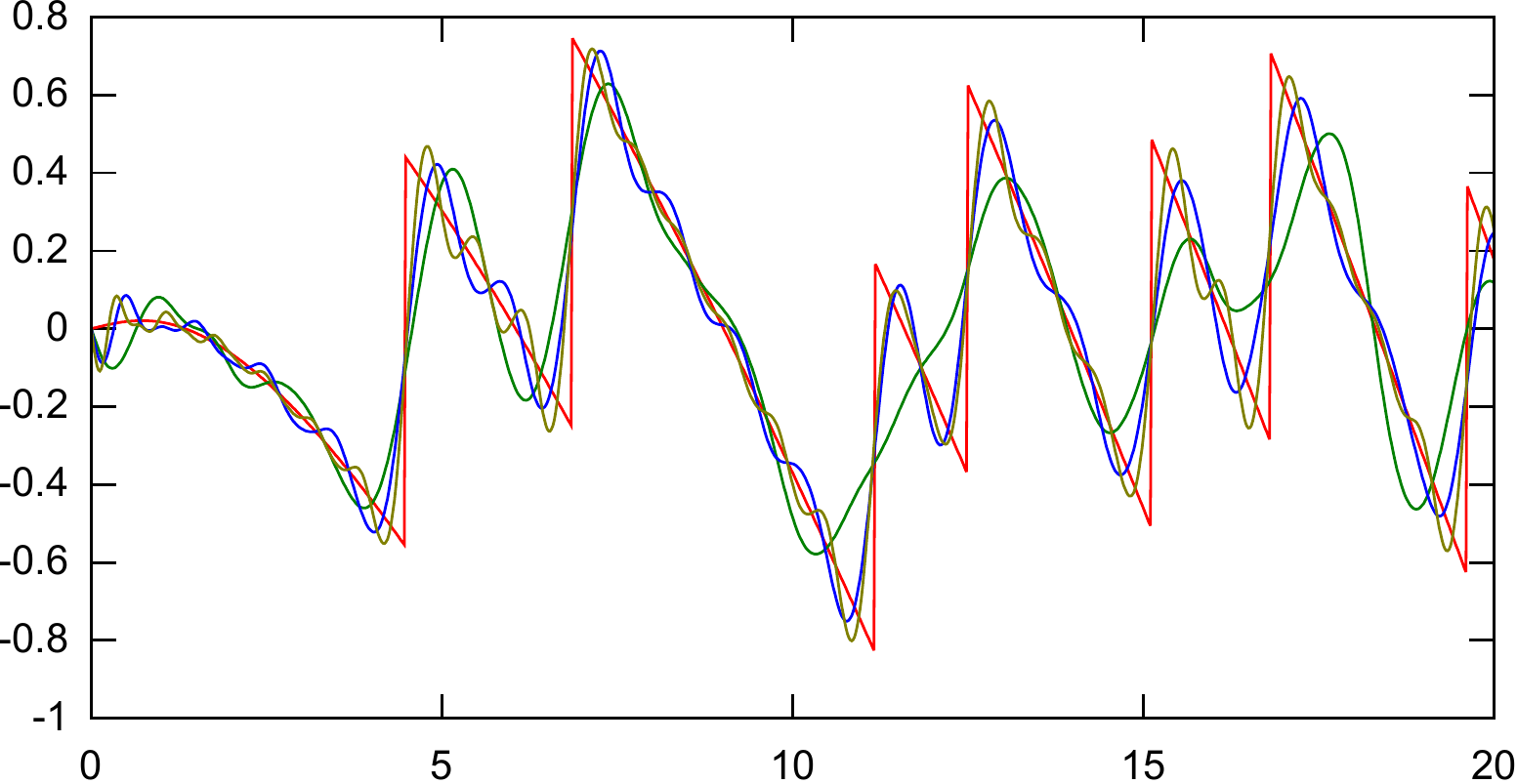}
  \end{center}
  \caption{$R_x(t)$ for $\chi_{7b}$}
  \label{Figure8}
 \end{minipage}
\end{figure}

Figures~\ref{Figure6}, \ref{Figure7} and \ref{Figure8} indicate the values
$$R_x(t)=\frac1\pi\Im\log L_{x}\l(\frac12+it,\,\chi\r)$$
for $\chi_{3}$, $\chi_{7a}$, $\chi_{7b}$, respectively,
for $x=p_{10}$ (green), $x=p_{100}$ (blue), $x=p_{1000}$ (yellow) and
$\infty$ (red).
This also seems to reflect the property of DRH.
The green, blue and yellow curves appear to converge to the red one
more smoothly only when $\chi^2\ne\1$ (Figure \ref{Figure7}).
In the other two cases, the curves oscillate many times near the origin.

The leaps in the red curves correspond to the zeros of $L(s,\chi)$.
We normalize that the jumps at zeros are equal to one.
This reflects the 
\remnum{conjecture}{41} that the multiplicity of such zeros 
\remnum{should be}{41} all one.
In other words, if we express their derivatives by the Dirac delta
function, the coefficients are one.

\begin{figure}[h]
 \begin{minipage}[c]{0.5\textwidth}
  \begin{center}
     \includegraphics[width=20em]{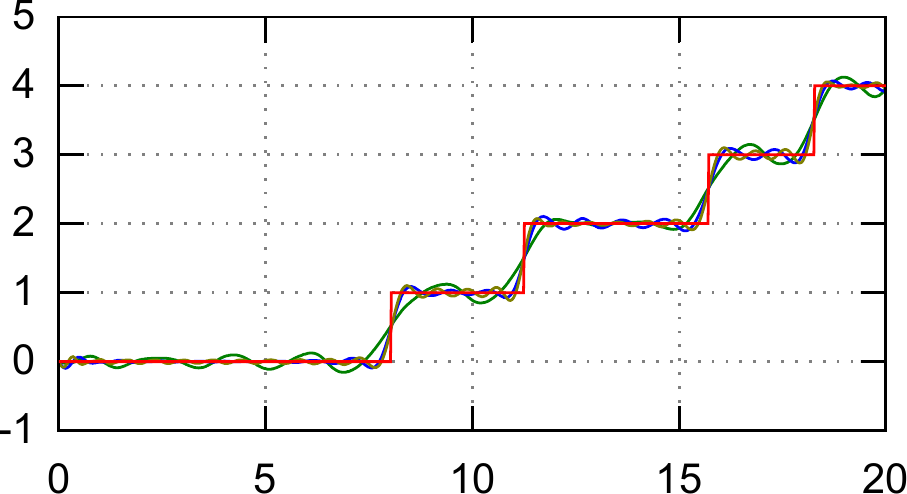}
  \end{center}
  \caption{$N_x(t)$ for $\chi_3$}
  \label{Num_func_3}
 \end{minipage}

 \vspace{1em}

 \begin{minipage}[c]{0.5\textwidth}
  \begin{center}
     \includegraphics[width=20em]{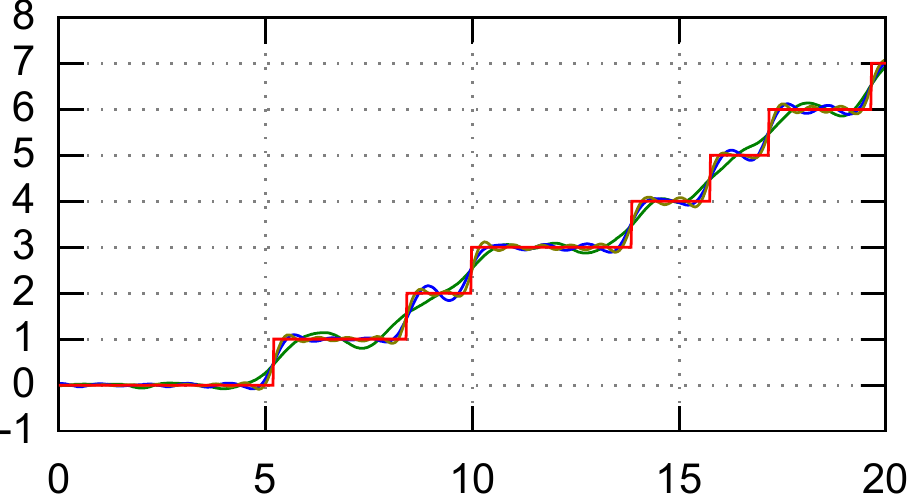}
  \end{center}
  \caption{$N_x(t)$ for $\chi_{7a}$}
  \label{Num_func_7a}
 \end{minipage}%
 \begin{minipage}[c]{0.5\textwidth}
  \begin{center}
     \includegraphics[width=20em]{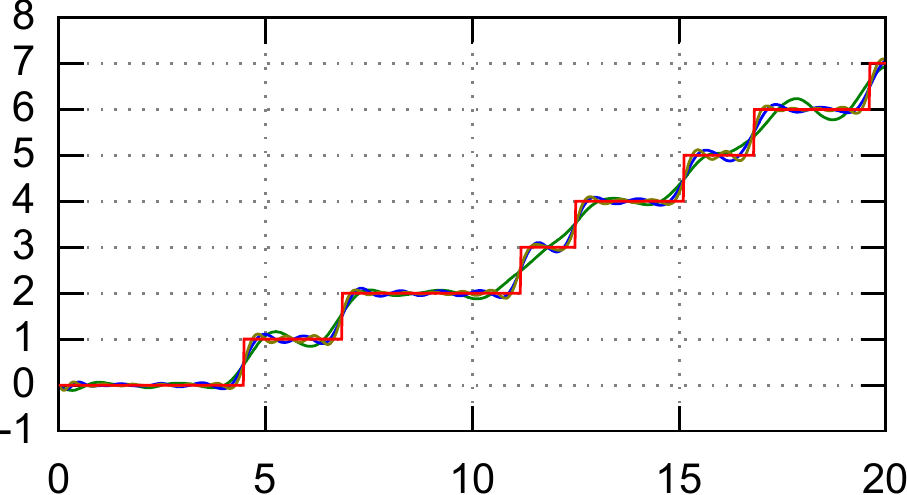}
  \end{center}
  \caption{$N_x(t)$ for $\chi_{7b}$}
  \label{Num_func_7b}
 \end{minipage}
\end{figure}

\begin{figure}[h]
 \begin{center}
  \includegraphics[width=18em]{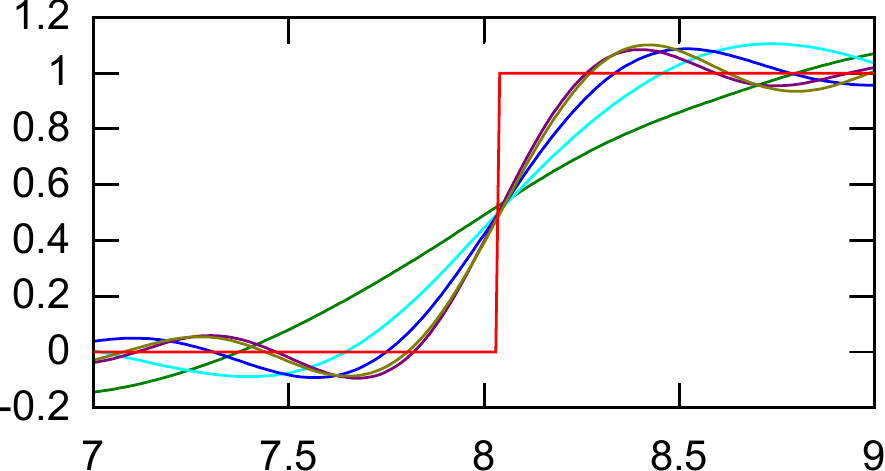} \qquad
  \includegraphics[width=18em]{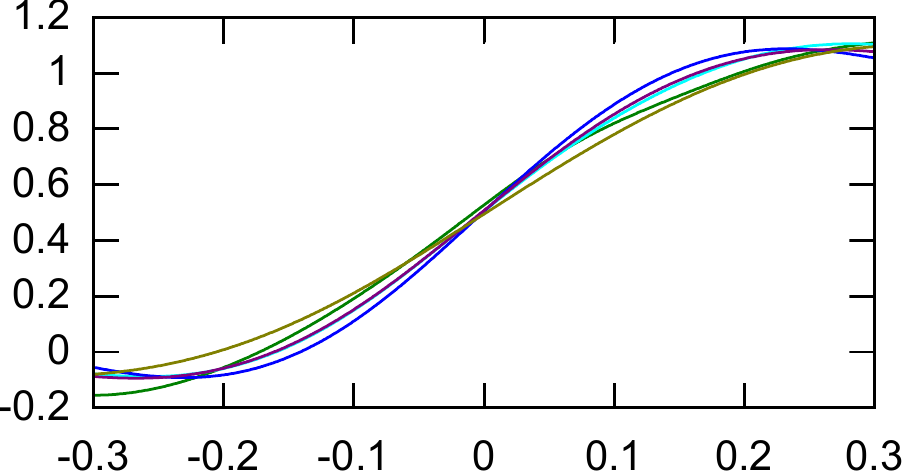}
 \end{center}
 \caption{$N_x(t)$ (left) and $\tilde{N}_x(z)$ (right) for $\chi_3$}
 \label{scale_mod3}
\end{figure}

\begin{figure}[h]
 \begin{center}
  \includegraphics[width=18em]{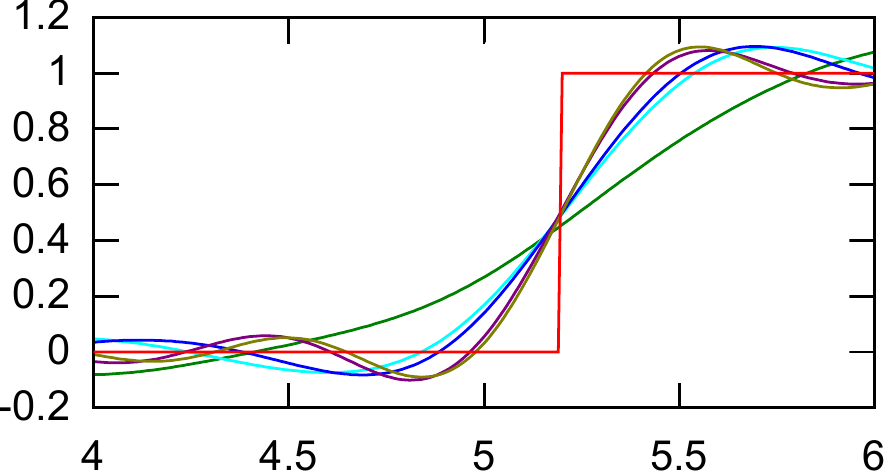} \qquad
  \includegraphics[width=18em]{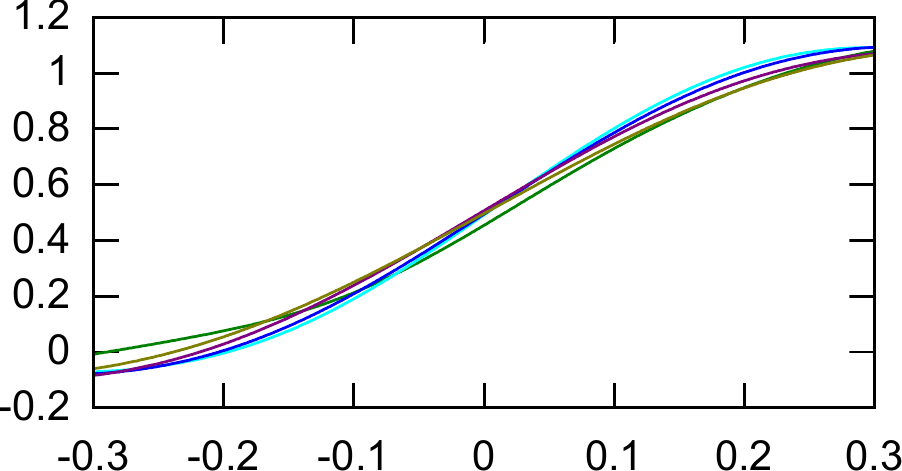}
 \end{center}
 \caption{$N_x(t)$ (left) and $\tilde{N}_x(z)$ (right) for $\chi_{7a}$}
 \label{scale_mod7a}
\end{figure}

\begin{figure}[h]
 \begin{center}
  \includegraphics[width=18em]{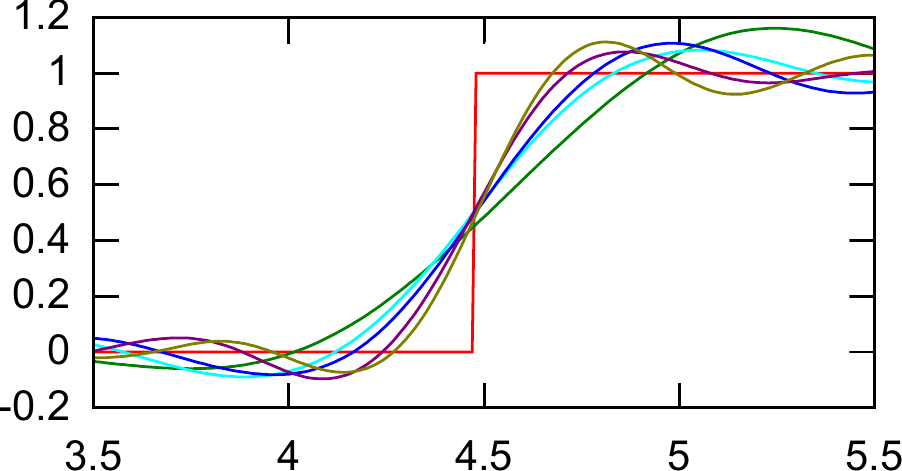} \qquad
  \includegraphics[width=18em]{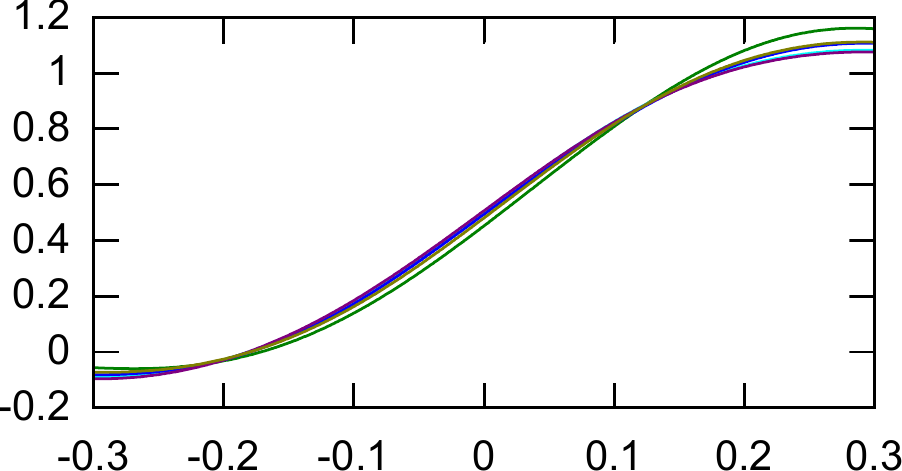}
 \end{center}
 \caption{$N_x(t)$ (left) and $\tilde{N}_x(z)$ (right) for $\chi_{7b}$}
 \label{scale_mod7b}
\end{figure}

We define another function $N_x(t)$ from $R_x(t)$ by
subtracting the contribution of the $L$-function versions of the
Riemann-Siegel theta function.
This counts the number of the nontrivial zeros on the critical line in
the limit of $x \to \infty$.
Figures~\ref{Num_func_3}, \ref{Num_func_7a} and \ref{Num_func_7b} show
the values of $N_x(t)$ for $\chi_3$, $\chi_{7a}$, $\chi_{7b}$,
respectively.
\remnum{The panels}{42} of Figures~\ref{scale_mod3}, \ref{scale_mod7a},
\ref{scale_mod7b} show $N_x(t)$ around the smallest nontrivial zeros of the
$L$-functions with $x = p_{10}$ (green), $x = p_{50}$ (light blue), $x =
p_{100}$ (blue), $x = p_{500}$ (purple), $x = p_{1000}$ (yellow) and $x
= \infty$ (red).
As the case of $R_x(t)$, we see a sharp step structure as the cut-off
parameter $x$ getting larger.

These figures also tell us that the values $\Im\log L(\frac12+it)$ are
almost stable for nontrivial zeros $\frac12+it$ of the $L$-function, 
no matter how many prime numbers we take into account.
This suggests that the nontrivial zeros are analogs of the critical
points in statistical mechanics, which are stable to the finite-size
correction.

\begin{table}[t]
 \begin{center}
  \begin{tabular}{ccccc} \hline \hline
   character && $t_1$ && $\lambda$ \\ \hline
   $\chi_3$ && 8.0397... && 0.217\\ 
   $\chi_{7a}$ && 5.1981... && 0.193 \\ 
   $\chi_{7b}$ && 4.4757... && 0.151 \\ 
   \hline \hline
  \end{tabular}
 \end{center}
  \caption{Numerically evaluated exponents around the smallest zeros
 $\frac{1}{2}+it_1$ for $\chi_3$, $\chi_{7a}$ and $\chi_{7b}$}
  \label{exponents}
\end{table}

To examine the analogy to critical phenomena in \remnum{statistical}{43}
mechanics, we shall check the scaling property around the critical point.
Being the smallest zero $\frac{1}{2}+it_1$, we define the scaling variable 
\[
  z = \frac{t - t_1}{t_1} x^\lambda \, .
\]
Correspondingly we introduce a {\em scaled} function $\tilde{N}_x(z)$,
defined as $N_x(t) = \tilde{N}_x(z = \frac{t - t_1}{t_1} x^\lambda)$.
Right panels of Figures~\ref{scale_mod3}, \ref{scale_mod7a} and
\ref{scale_mod7b} show the values of $\tilde{N}_x(z)$.
By choosing a proper exponent $\lambda$, all the curves are almost
approximated by only one curve.
This means that the dependence on the cut-off parameter $x$ appears only
in the form of the scaling variable $z$.
This scaling behavior supports the similarity to the critical phenomena.

Table~\ref{exponents} shows the numerical values of the smallest zeros
of the $L$-functions and the corresponding exponents for
$\chi_3$, $\chi_{7a}$ and $\chi_{7b}$.
These exponents are numerically determined by fitting the curves of
$\tilde{N}_x(z)$ by changing the parameter $x=10, 50, 100, 500, 1000$.

In the case of the ordinary critical phenomena, there is only one
critical point.
On the other hand, there are \remnum{infinitely many zeros}{44} on the
critical line of the $L$-function, which are analogs of the critical point.
Thus, even if we focus on only the smallest zero, as discussed in this
study, \remnum{there should be correction to its scaling bahavior from such
other zeros: we have to take care of the scaling property for others
simultaneously.}{45}

\end{document}